\documentclass{article}

\usepackage[round]{natbib}
\usepackage{arxiv}
\usepackage[T1]{fontenc}
\usepackage[english]{babel}
\usepackage{microtype}
\usepackage{indentfirst}
\usepackage{amsmath,amssymb,amsfonts,amsthm}
\usepackage{mathtools}
\usepackage{mathrsfs}
\usepackage{bm}
\usepackage{bbm}
\usepackage{enumitem}
\usepackage{textcomp}
\usepackage{manyfoot}
\usepackage{soul}
\usepackage{pifont}
\usepackage{fancybox}
\usepackage[title]{appendix}
\usepackage{graphicx}
\usepackage{caption}
\usepackage{subcaption}
\usepackage{xcolor}
\usepackage{booktabs}
\usepackage{tabularx}
\usepackage{array}
\usepackage{makecell}
\usepackage{multirow,multicol}
\usepackage{adjustbox}
\usepackage{algorithm,algorithmicx}
\usepackage{algpseudocode}
\usepackage{listings}
\usepackage{url}
\usepackage[colorlinks=true, allcolors=blue]{hyperref}
\usepackage[capitalize]{cleveref}

\theoremstyle{plain}
\newtheorem{theorem}{Theorem}[section]
\newtheorem{lemma}[theorem]{Lemma}
\newtheorem{corollary}[theorem]{Corollary}
\newtheorem{proposition}[theorem]{Proposition}

\theoremstyle{definition}
\newtheorem{definition}[theorem]{Definition}
\newtheorem{assumption}[theorem]{Assumption}

\theoremstyle{remark}

\Crefname{theorem}{Theorem}{Theorems}
\Crefname{lemma}{Lemma}{Lemmas}
\Crefname{corollary}{Corollary}{Corollaries}
\Crefname{proposition}{Proposition}{Propositions}
\Crefname{conjecture}{Conjecture}{Conjectures}
\Crefname{assumption}{Assumption}{Assumptions}
\Crefname{definition}{Definition}{Definitions}
\Crefname{example}{Example}{Examples}
\Crefname{problem}{Problem}{Problems}
\Crefname{exercise}{Exercise}{Exercises}
\Crefname{remark}{Remark}{Remarks}
\Crefname{figure}{Figure}{Figures}

\crefformat{equation}{(#2#1#3)}
\Crefformat{equation}{(#2#1#3)}
\crefrangeformat{equation}{(#3#1#4)--(#5#2#6)}
\Crefrangeformat{equation}{(#3#1#4)--(#5#2#6)}
\crefmultiformat{equation}{(#2#1#3)}{ and~(#2#1#3)}{, (#2#1#3)}{, and~(#2#1#3)}
\Crefmultiformat{equation}{(#2#1#3)}{ and~(#2#1#3)}{, (#2#1#3)}{, and~(#2#1#3)}

\graphicspath{{./figures/}}

\title
{Optimal strong convergence rates of explicit stabilized method for the Galerkin finite element discretization SPDEs driven by multiplicative noise}
\author{
	Junwei Huang\textsuperscript{1},
	Aiguo Xiao\textsuperscript{2},
	Xiao Tang\textsuperscript{1,*}\\[8pt]
	\begin{minipage}{0.9\textwidth}
		\centering\small
		\textsuperscript{1}School of Mathematics and Computational Science \&
		Hunan Key Laboratory for Computation and Simulation in Science and Engineering,
		Xiangtan University, Xiangtan 411105, China\\[4pt]
		\textsuperscript{2}Hunan Research Center of the Basic Discipline
		``Fundamental Algorithmic Theory and Novel Computational Methods'' \&
		China National Center for Applied Mathematics in Hunan,
		Xiangtan University, Xiangtan 411105, China\\[4pt]
		\textsuperscript{*}Corresponding author: tangx@xtu.edu.cn
	\end{minipage}
}
\date{}

\begin{document}

\maketitle

\begin{abstract}
Explicit stabilized methods provide an efficient approach to stiff systems arising from deterministic parabolic equations, yet rigorous convergence theory for
their application to stochastic partial differential
equations remains scarce.
To the best of our knowledge, the only prior fully discrete strong
convergence analysis of explicit stabilized integrators for SPDEs is
due to [A. Abdulle et al. IMA J. Numer. Anal. 43(1)(2023), pp. 258--292] and it is restricted to additive noise.
The present work addresses the extension to multiplicative noise, by coupling Galerkin finite element spatial discretization with explicit stabilized time integrators.
Under global Lipschitz assumptions and appropriate spatial
regularity conditions, we derive the strong error bound
$\mathcal{O}(h^{1+r}+\tau^{1/2})$ with $r\in[0,1)$,
which achieves the optimal temporal strong convergence rate
$1/2$ for the considered family of methods.
One- and two-dimensional numerical tests with multiplicative
noise validate the predicted convergence rates and compare
the computational performance of explicit stabilized methods
against the linearly implicit Euler method.
\end{abstract}
\keywords{Semilinear parabolic SPDEs\and Multiplicative noise\and Explicit stabilized methods\and Finite element methods\and Strong convergence.}

\section{Introduction}
Semilinear parabolic stochastic partial differential equations are important models describing the combined effects of diffusion, nonlinear reaction, and stochastic perturbations. We consider equations of the form \cite{DaPratoZabczyk2014}
\begin{equation}\label{eq:1.1}
	\begin{cases}
		\mathrm{d}u(t)+Au(t)\,\mathrm{d}t
		=
		F(u(t))\,\mathrm{d}t
		+
		G(u(t))\,\mathrm{d}W^Q(t),
		& t\in(0,T],\\
		u(0)=u_0,
	\end{cases}
\end{equation}
where \(H=L^2(D)\), with \(D\subset\mathbb{R}^d\) a bounded convex polygonal
(polyhedral) domain, and \(A:\mathcal{D}(A)\subset H\to H\) is a positive
definite, self-adjoint elliptic operator with homogeneous Dirichlet
boundary conditions. The operator \(-A\) generates an analytic semigroup
\[
E(t)=e^{-tA},\qquad t\ge 0.
\]
Moreover, $F$ denotes the nonlinear drift coefficient,
$G$ is the diffusion coefficient, and
$\{W^Q(t)\}_{t\ge0}$ is a possibly cylindrical
$Q$-Wiener process.
The precise assumptions on the coefficients
and initial value are stated in Section~\ref{sec2}.
Under the assumptions stated in Section~\ref{sec2}, equation \eqref{eq:1.1}
admits a unique mild solution satisfying
\begin{equation}\label{eq:1.2}
	u(t)
	=
	E(t)u_0
	+
	\int_0^t E(t-s)F(u(s))\,\mathrm{d}s
	+
	\int_0^t E(t-s)G(u(s))\,\mathrm{d}W^Q(s),
	\qquad t\in[0,T].
\end{equation}

A Galerkin finite element discretization of \eqref{eq:1.1} leads to a system of
stochastic differential equations.
For standard piecewise linear elements on quasi-uniform meshes of size $h$, the largest eigenvalue of the discrete elliptic operator $A_h$ grows like $h^{-2}$.
Consequently, the standard explicit Euler method requires $\tau=O(h^2)$ for diffusion stability,
which substantially limits its efficiency on fine meshes. To relax this restriction, implicit and linearly implicit methods are widely used for the numerical solution of parabolic stochastic partial differential equations.
In particular, the linearly implicit Euler method treats the diffusion term implicitly
by means of the resolvent operator $(I+\tau A_h)^{-1}$,
and is therefore free from the above explicit step-size constraint.
Nevertheless, it requires solving a linear system involving $A_h$ at every time step,
which incurs additional computational and storage costs
in high dimensions, and when simulating a large number of sample paths.

Explicit stabilized methods offer an alternative
to implicit methods for handling diffusion stiffness;
see \cite{Abdulle2015} for a review.
Chebyshev-based Runge--Kutta methods enlarge the stability interval along the
negative real axis by increasing the number of internal stages; a
non-exhaustive list of references is
\cite{SommeijerShampineVerwer1998,
	AbdulleMedovikov2001,Abdulle2002,TangXiao2020}.
Further developments include multiscale and
partitioned formulations
\cite{AbdulleVilmart2012,AbdulleVilmart2013,TangHuang2026,TangXiao2019},
adaptive strategies \cite{HuangTang2026},
and applications to optimization \cite{EftekhariEtAl2021}.
The stabilization idea has also been extended to stochastic differential
equations; see, e.g.,
\citealp{AbdulleCirilli2008,AbdulleLi2008,AbdulleVilmartZygalakis2013,YuXiaoTang2025,TangXiao2019},
among others.
In particular, 
\cite{AbdulleAlmuslimaniVilmart2018} introduced
the SK-ROCK method with favorable mean-square
stability properties.

The theory and numerical approximation of SPDEs
have been extensively studied; see
\cite{Printems2001,Kruse2014b,LordPowellShardlow2014,
	BrehierGoudenege2019,BrehierVilmart2016,
	DavieGaines2001,GyongyMillet2005,GyongyNualart1997,
	Walsh2005,Debussche2011,Hausenblas2003,
	JentzenKloeden2009OrderBarrier,JentzenKloedenWinkel2011,
	KovacsLarssonLindgren2010,Kruse2014Milstein,
	QiWangYan2025,Yan2005,DebusschePrintems2009,
	GeissertKovacsLarsson2009,PrevotRockner2007,
	Brehier2014,Brehier2024TotalVariation,
	CuiHong2019,BrehierCuiHong2019, BauzetSultanValletZimmermann2024,YangZhao2021,MukamTambue2020,FengLiZhang2021,BanasVieth2023}
for a non-exhaustive list of references.
For problems with multiplicative noise,
\cite{Kruse2014a} established optimal strong
error estimates for Galerkin finite element
approximations combined with the linearly implicit
Euler method.
\cite{LordTambue2013} analyzed
the strong convergence of stochastic exponential
integrators for finite element discretizations
with both additive and multiplicative noise.
For additive-noise problems,
\cite{Wang2017} derived strong error estimates
for the finite element discretization combined
with the linearly implicit Euler method,
achieving second-order spatial and first-order
temporal convergence under suitable regularity
assumptions.

For explicit stabilized methods,  \cite{ABV2023}
established a strong convergence framework
for finite element approximations with additive noise.
Their error bounds are uniform in the spatial
mesh size, the time step, and the number of
internal stages under a suitable stability condition.
Theorem~3.7 of \citet{ABV2023} yields the temporal strong rate
\(1-\epsilon\) for the linear problem \(F=0\) with spatially regular
additive noise, and Remark~4.7 of the same work sketches an extension to
a nonlinear drift \(F\) of class \(C^2\) with bounded first and second
derivatives. These results, however, are restricted to additive noise; for
multiplicative noise, no strong convergence theorem is established in
\citet{ABV2023}.

This paper derives the optimal strong convergence order for explicit stabilized finite element methods with multiplicative noise. Because the two operator functions \(\mathcal{A}_k(-\tau A_h)\) and
\(\mathcal{B}_k(-\tau A_h)\) of the explicit stabilized scheme
(defined in Section~\ref{sec3}) are polynomials in the discrete diffusion
operator \(A_h\), whereas the linearly implicit Euler scheme uses the
resolvent \((I+\tau A_h)^{-1}\), these polynomial propagators have
weaker smoothing properties, which introduces several challenges in
our analysis.  Such weaker regularization effects were also noted in \citet{ABV2023}, yet their analysis is restricted to additive noise. We  adopt a distinct analytical framework for the multiplicative noise setting. We establish pointwise and time-integrated error estimates for the discrete convolution operator, which hold uniformly with respect to the number of internal stages. Using these estimates, we are able to control the solution-dependent noise and recover the full spatial  order.

Under the same assumptions on the coefficients
and initial value as in \cite{Kruse2014a},
together with suitable conditions on the stabilized
propagators, we prove that, for $r\in[0,1)$
and $p\in[2,\infty)$,
\begin{equation}\label{eq:1.3}
	\|u(t_n)-u_h^n\|_{L^p(\Omega;H)}
	\le C\bigl(h^{1+r}+\tau^{1/2}\bigr),
	\qquad t_n=n\tau\le T.
\end{equation}
The constant $C$ is independent of $h$, $\tau$,
and the number of internal stages satisfying
the stability condition.
Thus, the explicit stabilized methods attain the same convergence orders as the linearly implicit Euler method, without strengthening the assumptions on the coefficients or initial value.

The remainder of this paper is organized as follows.
Section~\ref{sec2} introduces the analytical framework
and states the assumptions used throughout the paper.
Section~\ref{sec3} establishes pointwise and time-integrated
propagator error estimates uniform in the number
of internal stages.
Section~\ref{sec4} proves optimal strong convergence rates
for multiplicative noise.
Section~\ref{sec5} considers additive noise as a special case
and derives improved convergence rates under
additional assumptions, including first-order temporal
strong convergence without an $\varepsilon$-loss.
Finally, Section~\ref{sec6}  presents numerical experiments
with multiplicative noise to illustrate the theoretical
rates and compare the computational efficiency
with that of the linearly implicit Euler method.

\section{Setting and preliminaries}\label{sec2}

Let $H$ be the separable Hilbert space $L^2(\mathcal{D})$, where $\mathcal{D} \subset \mathbb{R}^d$ is  a bounded, open, and convex polygonal/polyhedral domain in dimension $d \in \{1,2,3\}$. The inner product in $H$ is denoted by $(\cdot, \cdot)$. The norm in $H$ is denoted by $\|\cdot\|$. We denote by $\mathcal{L}(H)$ the space of bounded linear operators on $H$, equipped with the operator norm $\|\cdot\|_{\mathcal{L}(H)}$,
and by $\mathcal{L}_2(H)$ the space of Hilbert--Schmidt operators on $H$,
equipped with the norm $\|\cdot\|_{\mathcal{L}_2(H)}$.

The assumptions stated in this section are standard in the literature on SPDEs \cite{DaPratoZabczyk2014} and their numerical approximations \cite{Kruse2014a,LordPowellShardlow2014}.

\subsection{Linear operator $A$}
Let $A : \mathcal{D}(A) \subset H \to H$ be a linear, densely defined, positive self-adjoint unbounded operator with compact inverse. Hence,  there exists a nondecreasing sequence of positive eigenvalues $(\lambda_n)_{n\geq1}$ and an orthonormal basis of eigenvectors $(e_{n})_{n\geq1}$ in $H$ such that $Ae_n = \lambda_ne_n$ and
$$0 < \lambda_1 \leq \lambda_2 \leq \dots \leq \lambda_n~(\to \infty).$$
The domain of $A$ is characterized by
$$
\mathcal D(A) = \left\{ x \in H : \sum_{n=1}^\infty \lambda_n^2 (x, e_n)^2 < \infty \right\}.
$$
Consequently, $-A$ generates an analytic contraction semigroup $(E(t))_{t\geq0}$ on $H$, given by $ E(t)=e^{-tA}, t\geq0.$
For the purposes of the subsequent analysis, we introduce the notion of fractional powers of the operator $A$, which is formulated below. For any $r \in \mathbb{R}$, the operator $A^{\frac{r}{2}} : D(A^{\frac{r}{2}}) \to H$ is given by
$$
A^{\frac{r}{2}} x = \sum_{n=1}^\infty \lambda_n^{\frac{r}{2}} x_n e_n
$$
for all
$$
x \in D\left(A^{\frac{r}{2}}\right) = \left\{ x = \sum_{n=1}^\infty x_n e_n : (x_n)_{n \ge 1} \subset \mathbb{R} \text{ with } \|x\|_r^2 := \sum_{n=1}^\infty \lambda_n^r x_n^2 < \infty \right\}.
$$
By defining $\dot{H}^r := D\left(A^{\frac{r}{2}}\right)$ together with the norm $\|\cdot\|_r$ for $r \in \mathbb{R}$, $\dot{H}^r$ becomes a Hilbert space. Note that we have $\|x\|_r = \left\|A^{\frac{r}{2}} x\right\|$ for all $r \in \mathbb{R}$ and $x \in \dot{H}^r$.

We next recall several standard regularity estimates for the analytic
semigroup $(E(t))_{t\geq0}$; see, for instance, \cite{Kruse2014a}.

\begin{proposition}
	\label{pro1}
	The analytic semigroup $(E(t))_{t\geq0}$ satisfies the following estimates.
	Every $\mu\geq0$,
	\begin{equation}
		\|A^\mu E(t)\|
		\leq Ct^{-\mu},
		\qquad t>0.
		\label{eq2.1}
	\end{equation}
	For every $\nu\in[0,1]$,
	\begin{equation}
		\|A^{-\nu}(I-E(t))\|
		\leq Ct^\nu,
		\qquad t>0.
		\label{eq2.2}
	\end{equation}
	Moreover, for every $x\in H$, $0\leq\tau_1<\tau_2$, and
	$\nu\in[0,1]$,
	\begin{equation}
		\int_{\tau_1}^{\tau_2}
		\left\|
		A^{\nu/2}E(\tau_2-\sigma)x
		\right\|^2\,d\sigma
		\leq
		C(\tau_2-\tau_1)^{1-\nu}\|x\|^2,
		\label{eq2.3}
	\end{equation}
	and
	\begin{equation}
		\left\|
		A^\nu
		\int_{\tau_1}^{\tau_2}
		E(\tau_2-\sigma)x\,d\sigma
		\right\|
		\leq
		C(\tau_2-\tau_1)^{1-\nu}\|x\|.
		\label{eq2.4}
	\end{equation}
\end{proposition}

Throughout the paper, $C$ denotes a generic positive constant,
which may vary from line to line and is independent of the discretization
parameters.

\subsection{Nonlinear coefficients $F$}

Fix a regularity parameter $r\in[0,1)$.
Following \cite{Kruse2014a}, we impose the following
assumptions on the drift and diffusion coefficients.

\begin{assumption}
	\label{non}
	The nonlinear operator $F:H\to\dot{H}^{r-1}$
	is globally Lipschitz continuous. More precisely,
	there exists a constant $L>0$ such that
	\begin{equation}\label{eq:F-Lipschitz}
		\|F(u)-F(v)\|_{r-1}
		\le L\|u-v\|,
		\qquad u,v\in H.
	\end{equation}
\end{assumption}
This assumption is weaker than global Lipschitz
continuity from $H$ to $H$, since drift increments
are measured in the negative-order norm
$\|\cdot\|_{r-1}$.
It also implies the linear growth bound
$\|F(v)\|_{r-1}\le C(1+\|v\|)$ for all $v\in H$.
No differentiability of $F$ is required here;
additional assumptions for the enhanced additive-noise
analysis are introduced in Section~5.

\subsection{$Q$-Wiener process and diffusion coefficient}
Before stating the assumptions on the noise term, we first introduce the $Q$-Wiener process $(W^{Q}(t))_{t\geq0}$, which will be used throughout the rest of the paper. Let $(\epsilon_n)_{n\ge 1}$ be a complete orthonormal system of the Hilbert space $H$, and consider a bounded sequence $(q_n)_{n\ge 1}$ of non-negative real numbers. Let $Q$ and $Q^\frac{1}{2}$ be the bounded, linear, self-adjoint operators on $H$, defined by
$$
Qu = \sum_{n=1}^\infty  q_n ( u, \epsilon_n ) \epsilon_n,
\quad
Q^{\frac{1}{2}} u = \sum_{n=1}^\infty  \sqrt{q_n} ( u, \epsilon_n ) \epsilon_n.
$$

\begin{definition}
	Let
	$(\Omega,\mathcal{F},(\mathcal{F}_t)_{t\geq0},
	\mathbb{P})$
	be a filtered probability space satisfying
	the usual conditions.
	Let $(\beta_n)_{n\geq1}$ be a family of independent
	standard real-valued $(\mathcal{F}_t)$-Wiener processes.
	The possibly cylindrical $Q$-Wiener process
	is represented by the formal series
	$$
	W^Q(t)
	=
	\sum_{n=1}^{\infty}
	\sqrt{q_n}\,\beta_n(t)\epsilon_n,
	\qquad t\geq0.
	$$
\end{definition}
If $\operatorname{Tr}(Q)<\infty$, the above series converges
in $L^2(\Omega;H)$ for each fixed $t\geq 0$ and defines an
$H$-valued $Q$-Wiener process. Otherwise, it is interpreted
as a cylindrical Wiener process with covariance operator $Q$.

For the stochastic integration framework, we introduce
$H_0:=Q^{1/2}(H)$ with the inner product
\[
(u,v)_{H_0}
:= (Q^{-1/2}u,Q^{-1/2}v),
\]
where $Q^{-1/2}$ denotes the pseudoinverse of $Q^{1/2}$.
Let $\mathcal{L}_2^0:=\mathcal{L}_2(H_0,H)$ denote the space
of Hilbert--Schmidt operators from $H_0$ into $H$, equipped
with the norm
$$
\|\Phi\|_{\mathcal{L}_2^0}^2
= \|\Phi Q^{1/2}\|_{\mathrm{HS}}^2
= \sum_{i=1}^{\infty}\|\Phi\psi_i\|^2,
$$
where $(\psi_i)_{i\geq1}$ is any orthonormal basis of $H_0$.

Using the fractional powers of $A$ introduced above,
we impose the following Lipschitz and spatial
regularity conditions on the diffusion coefficient $G$.
\begin{assumption}
	\label{ass:noise}
	Let $\{W^Q(t)\}_{t\in[0,T]}$ be the possibly
	cylindrical $Q$-Wiener process introduced above,
	and let $r\in[0,1)$ be the regularity parameter
	in Assumption~\ref{non}.
	The nonlinear operator $G:H\to\mathcal{L}_2^0$
	satisfies
	\begin{equation}\label{eq:G-Lipschitz}
		\|G(u)-G(v)\|_{\mathcal{L}_2^0}
		\le L\|u-v\|,
		\qquad u,v\in H,
	\end{equation}
	for some constant $L>0$.
	Furthermore, for every $v\in\dot{H}^{r}$,
	$A^{r/2}G(v)\in\mathcal{L}_2^0$ and
	\begin{equation}\label{eq:G-regularity}
		\|A^{r/2}G(v)\|_{\mathcal{L}_2^0}
		\le L(1+\|v\|_r),
		\qquad v\in\dot{H}^{r}.
	\end{equation}
\end{assumption}
No differentiability assumptions on $G$ are imposed.
Additional assumptions used to obtain improved
temporal convergence rates for additive noise
are stated in Section~5.

\subsection{Initial value and solution regularity}
\begin{assumption}
	\label{ass:init}
	For a given \(p\in[2,\infty)\), the initial value
	\(u_0:\Omega\to\dot H^{1+r}\) is an \(\mathcal{F}_0\)-measurable
	random variable satisfying \(u_0\in L^p(\Omega;\dot H^{1+r})\).
\end{assumption}

We assume spatially regular initial data to simplify
the analysis and avoid singular factors near $t=0$
in the error estimates.
Less regular initial data can be considered using
semigroup smoothing and nonsmooth data error estimates,
but the resulting bounds generally contain factors
that become singular as $t\to0$.
The constants in our estimates may depend on $T$,
but are independent of the discretization parameters.

Under the assumptions stated above,
equation~\eqref{eq:1.1} admits a unique mild solution
with the following regularity properties.

\begin{theorem}
	\label{th1}
	Let $r\in[0,1)$ and $p\in[2,\infty)$.
	Suppose that the assumptions on $A$ and
	Assumptions~\ref{non}, \ref{ass:noise},
	and~\ref{ass:init} hold for these parameters.
	Then equation~\eqref{eq:1.1} admits a unique
	mild solution given by \eqref{eq:1.2}, satisfying
	\begin{equation}\label{eq2.10}
		\sup_{t\in[0,T]}
		\|u(t)\|_{L^p(\Omega;\dot{H}^{1+r})}
		\le C.
	\end{equation}
	Furthermore, for every $\sigma\in[0,1+r]$,
	there exists a constant $C_\sigma$ such that
	\begin{equation}\label{eq2.11}
		\|u(t)-u(s)\|_{L^p(\Omega;\dot{H}^{\sigma})}
		\le C_\sigma
		|t-s|^{\min\{1/2,(1+r-\sigma)/2\}},
		\qquad 0\le s<t\le T.
	\end{equation}
	Here $C$ and $C_\sigma$ are independent of the discretization parameters.
\end{theorem}
\begin{proof}
	Existence and uniqueness follow from the
	well-posedness theory recalled in
	\cite{KruseLarsson2012}.
	The spatial bound and the temporal estimates
	for $\sigma\in[0,1+r)$ follow from
	\cite{KruseLarsson2012}.
	The case $\sigma=1+r$ follows directly
	from the uniform bound \eqref{eq2.10}.
\end{proof}
Combining the growth conditions on $F$ and $G$
with Theorem~\ref{th1} yields the following bounds.
\begin{corollary}
	Under the assumptions of Theorem~\ref{th1},
	there exists a constant $C$ such that
	\begin{equation}\label{eq2.12}
		\sup_{t\in[0,T]}
		\left(
		\|F(u(t))\|_{L^p(\Omega;\dot{H}^{r-1})}
		+
		\|A^{r/2}G(u(t))\|_{L^p(\Omega;\mathcal{L}_2^0)}
		\right)
		\le C.
	\end{equation}
	The constant $C$ is independent
	of the discretization parameters.
\end{corollary}
\begin{proof}
	The assertion follows from the growth bounds
	on $F$ and $G$ and Theorem~\ref{th1}.
\end{proof}

\section{Numerical methods} \label{sec3}
This section introduces the semi-discrete and fully discrete approximations of the process $u(t)$. For spatial discretization, we employ the Galerkin finite element method, while an explicit stabilized scheme is used for the time discretization.

\subsection{Galerkin finite element method}
In this section, we first review the Galerkin finite element method used for the spatial discretization of the SPDE. Let the spatial domain $\mathcal{D}$ be a bounded convex polygonal or polyhedral domain. We introduce a quasi-uniform triangulation $\mathcal{T}_h$ of $\mathcal{D}$ with mesh size $h\in(0,1]$. The finite element space $V_h$ is a finite-dimensional subspace of $V=H_0^1(\mathcal{D})$,
with dimension $N_h$. For simplicity, we take $V_h$ to consist of continuous, piecewise affine functions on $\mathcal{T}_h$ that vanish on the boundary $\partial\mathcal{D}$.

For the finite element discretization considered here,
we take $A=-\Delta$ with homogeneous Dirichlet
boundary conditions on $\mathcal{D}$.
Then
$\dot{H}^1=H_0^1(\mathcal{D})$ and
$\dot{H}^2=H^2(\mathcal{D})\cap H_0^1(\mathcal{D})$,
with equivalent norms.
For $u_h \in V_h$ we define $A_h u_h$ to be the unique element in $V_h$ which satisfies the relationship
\begin{equation*}
	a(u_h, y_h) := (A^{\frac{1}{2}}u_h,A^{\frac{1}{2}}y_h) = (A_h u_h, y_h) \quad \forall y_h \in V_h \subset \dot{H}^1 ,
\end{equation*}
where \( a : \dot{H}^1 \times \dot{H}^1 \to \mathbb{R} \) is the bilinear form corresponding to \( A \).

Since
$$(A_hu_h,y_h)=a(u_h,y_h)=(u_h,A_hy_h), \quad \forall u_h,y_h\in V_h,$$
and 
$$ (A_hu_h,u_h)=a(u_h,u_h)=\|u_h\|_1^2>0, \quad  \forall u_h\in V_h\setminus\{0\},$$
the operator $A_h$ is self-adjoint and positive definite
on $V_h$ with respect to the inner product of $H$. Consequently, there exist an $H$-orthonormal basis $\{e_{h,n}\}_{n=1}^{N_h}$ of $V_h$ and positive eigenvalues $\{\lambda_{h,n}\}_{n=1}^{N_h}$ such that $A_he_{h,n}=\lambda_{h,n}e_{h,n}, 1\leq n\leq N_h, $ with $0<\lambda_{h,1}\leq\lambda_{h,2} \leq\cdots\leq\lambda_{h,N_h}.$ Hence, $-A_h$ generates an analytic contraction semigroup
on $V_h$, denoted by $E_h(t)=e^{-tA_h}$. For every $\rho\geq0$, the following smoothing estimate holds
\begin{equation}
	\|A_h^\rho E_h(t)y_h\|
	\leq C_\rho t^{-\rho}\|y_h\|,
	\qquad t>0,\quad y_h\in V_h,
	\label{eq:3.1.1}
\end{equation}
where $C_\rho$ is independent of $h$ and $t$.

Next, we define a generalized orthogonal projector, \( P_h : \dot{H}^{-1} \to V_h \):
$$
( P_h u, u_h) = \langle u, u_h \rangle \quad \forall u \in \dot{H}^{-1}, \quad u_h \in V_h,
$$
where $\langle \cdot, \cdot \rangle = a(A^{-1} \cdot, \cdot) $ denotes the duality pairing between $\dot{H}^{-1}$ and $\dot{H}^{1}$.  For the quasi-uniform finite element meshes considered
here, the projector $P_h$ satisfies
\begin{flalign}
	\|A_h^{1/2}P_hv\|
	\le C\|v\|_1,
	\qquad v\in\dot{H}^1. \label{eq3.1.1}
\end{flalign}
Moreover, its extension to $\dot{H}^{-1}$ satisfies
\begin{flalign}
	\|A_h^{-1/2}P_hv\|
	\le \|v\|_{-1},
	\qquad v\in\dot{H}^{-1}. \label{eq3.1.2}
\end{flalign}
Therefore, the standard finite element method for \eqref{eq:1.1} can be written as: for \( h \in (0, 1] \),
\begin{equation}\label{eq:3.1.2}
	\left\{
	\begin{aligned}
		\mathrm{d}u_h(t)
		&=
		\bigl(-A_hu_h(t)+P_hF(u_h(t))\bigr)\,\mathrm{d}t
		+P_hG(u_h(t))\,\mathrm{d}W^Q(t),\\[2pt]
		u_h(0)&=P_hu_0.
	\end{aligned}
	\right.
\end{equation}

\subsection{Explicit stabilized methods}
\label{sec3.2}
We describe the explicit stabilized time discretization and the two Chebyshev-based schemes considered below; further details can be found in \cite{ABV2023,AbdulleAlmuslimaniVilmart2018,EftekhariEtAl2021}.
Let $\tau=T/M$, $t_m=m\tau$, and
$\Delta W_m^Q=W^Q(t_{m+1})-W^Q(t_m)$.
For a given number of internal stages $k\geq1$, the fully
discrete scheme is
\begin{flalign} \label{eq:3.2.1}
	\left\{
	\begin{aligned}
		u_{h,m+1}
		&=
		\mathcal{A}_k(-\tau A_h)u_{h,m}
		+\tau\mathcal{B}_k(-\tau A_h)
		P_hF(u_{h,m})
		+\mathcal{B}_k(-\tau A_h)
		P_hG(u_{h,m})\Delta W_m^Q, \\ 
		u_{0,h} &= P_hu_0,
		\qquad m=0,1,\ldots,M-1.
	\end{aligned}
	\right.
\end{flalign}

The operator functions are defined by
\begin{equation}
	f(-\tau A_h)v_h
	=
	\sum_{n=1}^{N_h}
	f(-\tau\lambda_{h,n})(v_h,e_{h,n})e_{h,n},
	\qquad
	f\in\{\mathcal{A}_k,\mathcal{B}_k\},
	\quad v_h\in V_h.
	\label{eq:3.2.2}
\end{equation}
Iterating \eqref{eq:3.2.1} gives the discrete mild formulation
\begin{equation}\label{eq:3.2.3}
	\begin{aligned}
		u_{h,m}
		={}&
		\mathcal{A}_k(-\tau A_h)^mP_hu_0\\
		&+\tau\sum_{i=0}^{m-1}
		\mathcal{A}_k(-\tau A_h)^{m-1-i}
		\mathcal{B}_k(-\tau A_h)P_hF(u_{i,h})\\
		&+\sum_{i=0}^{m-1}
		\mathcal{A}_k(-\tau A_h)^{m-1-i}
		\mathcal{B}_k(-\tau A_h)
		P_hG(u_{i,h})\Delta W_i^Q.
	\end{aligned}
\end{equation}
for $m=1,2,\ldots,M$.
The stage number $k$ is chosen such that
\begin{equation}
	\tau\lambda_{h,N_h}\leq L_k,
	\qquad
	\sup_{z\in[-L_k,0]}|\mathcal{A}_k(z)|\leq1.
	\label{eq:3.2.4}
\end{equation}

To obtain a stability interval whose length grows
quadratically with $k$, while providing damping away from
the origin, we use
\begin{equation}
	\mathcal{A}_k(z)
	=
	\frac{T_k(\omega_0+\omega_1z)}{T_k(\omega_0)},
	\label{eq:3.2.6}
\end{equation}
where $T_k$ is the Chebyshev polynomial of the first kind and
\begin{equation}
	\omega_0=1+\frac{\widetilde{\eta}}{k^2},
	\qquad
	\omega_1=\frac{T_k(\omega_0)}{T_k'(\omega_0)}.
	\label{eq:3.2.7}
\end{equation}
Here $\widetilde{\eta}>0$ is a fixed, sufficiently small
damping parameter. We take $L_k=2/\omega_1$, which grows
quadratically with $k$; see \cite{ABV2023}.

\begin{figure}[htbp]
	\centering
	\includegraphics[width=0.9\textwidth]{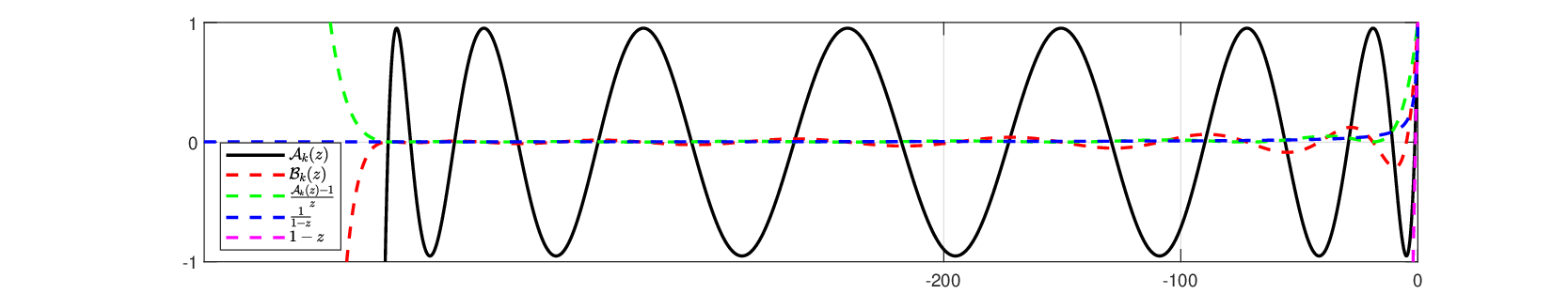}
	\caption{Stability functions $\mathcal{A}_k(z)$ and
		the two choices of $\mathcal{B}_k(z)$ for $k=15$
		and $\widetilde{\eta}=0.05$.}
	\label{fig:dampfunc}
\end{figure}

The first choice, corresponding to the SK-ROCK scheme
considered here, is
\begin{equation}
	\mathcal{B}_k(z)
	=
	\frac{U_{k-1}(\omega_0+\omega_1z)}{U_{k-1}(\omega_0)}
	\left(1+\frac{\omega_1z}{2}\right),
	\label{eq:3.2.8}
\end{equation}
where $U_{k-1}$ is the Chebyshev polynomial of the second kind.
For implementation, set
$$
H_{h,m}^{0}=u_{h,m},
\qquad
Z_{h,m}
=
\tau P_hF(u_{h,m})
+
P_hG(u_{h,m})\Delta W_m^Q.
$$
The scheme can then be evaluated through the recurrence \cite{AbdulleAlmuslimaniVilmart2018}
\begin{equation}
	\left\{
	\begin{aligned}
		H_{h,m}^{1}
		&=H_{h,m}^{0}
		-\tau\mu_1A_h
		\bigl(H_{h,m}^{0}+\nu_1Z_{h,m}\bigr)
		+\kappa_1Z_{h,m},\\
		H_{h,m}^{i}
		&=-\tau\mu_iA_hH_{h,m}^{i-1}
		+\nu_iH_{h,m}^{i-1}
		+\kappa_iH_{h,m}^{i-2},
		\qquad i=2,3,\ldots,k,\\
		u_{h,m+1}&=H_{h,m}^{k},
	\end{aligned}
	\right.
	\label{eq:3.2.9}
\end{equation}
where
$\mu_1=\omega_1/\omega_0$,
$\nu_1=k\omega_1/2$, and
$\kappa_1=k\mu_1$.
For $i=2,3,\ldots,k$, the coefficients are
\begin{equation}
	\mu_i=\frac{2\omega_1T_{i-1}(\omega_0)}{T_i(\omega_0)},
	\qquad
	\nu_i=\frac{2\omega_0T_{i-1}(\omega_0)}{T_i(\omega_0)},
	\qquad
	\kappa_i=-\frac{T_{i-2}(\omega_0)}{T_i(\omega_0)}
	=1-\nu_i.
	\label{eq:3.2.10}
\end{equation}
The second choice, referred to as the SK-ROCK variant, is
\begin{equation}
	\mathcal{B}_k(z)
	=
	\frac{\mathcal{A}_k(z)-1}{z},
	\qquad z\neq0,
	\qquad
	\mathcal{B}_k(0)=1.
	\label{eq:3.2.11}
\end{equation}
With the same initialization , its recurrence reads \cite{EftekhariEtAl2021}
\begin{equation}
	\left\{
	\begin{aligned}
		H_{h,m}^{1}
		&=H_{h,m}^{0}
		-\tau\mu_1A_hH_{h,m}^{0}
		+\mu_1Z_{h,m},\\
		H_{h,m}^{i}
		&=-\tau\mu_iA_hH_{h,m}^{i-1}
		+\mu_iZ_{h,m}
		+\nu_iH_{h,m}^{i-1}
		+\kappa_iH_{h,m}^{i-2},
		\qquad i=2,3,\ldots,k,\\
		u_{h,m+1}&=H_{h,m}^{k},
	\end{aligned}
	\right.
	\label{eq:3.2.12}
\end{equation}
where the remaining coefficients are given by
\eqref{eq:3.2.10}. These recurrences allow both schemes to be implemented
using only explicit stage computations. Fig.~\ref{fig:dampfunc} illustrates
$\mathcal{A}_k$ and the two choices of $\mathcal{B}_k$
along the negative real axis.

\subsection{Several assumptions and main result}

Since the subsequent convergence analysis is carried out under the following abstract conditions on the families  $(\mathcal{A}_{k})_{k\geq 1}$ and $(\mathcal{B}_{k})_{k\geq 1}$, we  give some assumptions on these families.
\begin{assumption}\label{ass:AB}
	For every \(k\ge 1\), let \(\mathcal{A}_k\) and \(\mathcal{B}_k\) be meromorphic functions. Suppose that there exist a sequence \((L_k)_{k\ge1}\) of positive real numbers, a constant \(C>0\) and a number \(\delta\in(0,1]\) such that
	\begin{align}
		&\mathcal{A}_k(0)=\mathcal{A}_k'(0)=\mathcal{B}_k(0)=1
		\quad\text{for all } k\ge1,\label{eq:3.3.1}\\
		&\sup_{k\ge1}\,\sup_{z\in[-L_k,0]}
		(1+\lvert z\rvert)\,\lvert\mathcal{B}_k(z)\rvert^{2}\le C,
		\label{eq:3.3.1.1}\\
		&\sup_{k\ge1}\,\sup_{z\in[-L_k,-\delta]}
		\lvert\mathcal{A}_k(z)\rvert<1,	\quad\text{for all } \delta \in (0,1]. \label{eq:3.3.2} 
	\end{align}
	Moreover,  there exists $\delta \in (0,1]$ such that
	\begin{align}
		&\sup_{k\ge1}\,\sup_{z\in\mathbb{D}(0,\delta)}
		\bigl(\lvert\mathcal{A}_k(z)\rvert+\lvert\mathcal{B}_k(z)\rvert\bigr)\le C, \quad\text{for all } \delta \in (0,1],
		\label{eq:3.3.3}
	\end{align}
	where \(\mathbb{D}(0,\delta)=\{z\in\mathbb{C}:\lvert z\rvert<\delta\}\) denotes the open disc of radius \(\delta\) centred at \(0\).
\end{assumption}
The two Chebyshev-based schemes introduced above
satisfy Assumption~\ref{ass:AB};
see \cite{ABV2023}.

%
	%

\begin{lemma}
	\label{lemma3.1}
	Let Assumption \ref{ass:AB} be satisfied. Then
	\begin{flalign}
		&(i)\sup_{k \geq 1, z \in [-L_k, 0]} \frac{|1 - \mathcal{B}_k(z)|}{|z|} \leq C,   \nonumber \\ 
		&(ii) \sup_{k \geq 1, z \in [-L_k, 0]} \frac{|e^{z} - \mathcal{B}_k(z)|}{|z|} \leq C, \nonumber \\
		&(iii) \sup_{k \geq 1, z \in [-L_k, 0]} \frac{|e^{z} - \mathcal{B}_k(z)|}{|z|^{\frac{1}{2}}} \leq C,     \nonumber \\
		&(iv)\sup_{j \geq 1, k \geq 1, z \in [-L_k, 0] } j 
		|\mathcal{A}_k(z)^j - e^{jz}| \leq C, \nonumber \\
		&(v)\sup_{j \geq 1, k \geq 1, z \in [-L_k, 0] } j^{1/2} |z|^{1/2} |\mathcal{A}_k(z)^{j-1} \mathcal{B}_k(z)| \leq C \label{eq:3.3.4}
	\end{flalign}
\end{lemma}
\begin{proof}
	The estimates \eqref{eq:3.3.4}$(i)$ and \eqref{eq:3.3.4}$(iv)$ can be found in \cite{ABV2023}. We prove the remaining three estimates. Let \( \widetilde{r}_k(z) = (e^{z} - \mathcal{B}_k(z))/z \). Since \( \mathcal{B}_k \) is a meromorphic function (by Assumption~\ref{ass:AB}), and $e^{0}-\mathcal{B}_k(0)=0$ by \eqref{eq:3.3.1}, while \( \mathcal{B}_k \) is bounded on \( \mathbb{D}(0, \delta) \) by \eqref{eq:3.3.3}, it follows that \( \widetilde{r}_k(z) \) is analytic in the open disc \( \mathbb{D}(0, \delta) \), where \( \delta \) is given by \eqref{eq:3.3.3} (recall that \( \delta \) does not depend on \( k \)). Owing to the maximum principle for analytic functions,
	$$
	\sup_{k\geq1,\ |z|\leq\delta/2}
	|\widetilde{r}_k(z)|
	\leq
	\frac{2}{\delta}
	\sup_{k\geq1,\ |z|=\delta/2}
	\bigl(e^{\delta/2}+|\mathcal{B}_k(z)|\bigr)
	\leq C.
	$$
	In addition, for $z \in [-L_{k},-\delta]$ in conjunction with \eqref{eq:3.3.3} we obtain
	$$
	\sup_{z \in [-L_k, -\delta/2]} |\widetilde{r}_k(z)| \leq \delta^{-1} \left( e^{-\delta} + \sup_{z \in [-L_k, -\delta/2]} |\mathcal{B}_k(z)| \right) \leq C.
	$$
	This completes the proof of \eqref{eq:3.3.4}$(ii)$. 
	Once \eqref{eq:3.3.4}$(iii)$ is established, \eqref{eq:3.3.4}$(ii)$ follows immediately by interpolation. Notice that for $z \in [-L_k, 0]$, $|e^z - \mathcal{B}_k(z)| \le 1 + \sup_{k \geq 1, z \in [-L_k, 0]} |\mathcal{B}_k(z)| \le C_1$. From \eqref{eq:3.3.4}$(ii)$, we also have $|e^z - \mathcal{B}_k(z)| \le C_2|z|$. Therefore,
	$$|e^z - \mathcal{B}_k(z)| = |e^z - \mathcal{B}_k(z)|^{\frac{1}{2}} \cdot |e^z - \mathcal{B}_k(z)|^{\frac{1}{2}} \le (C_2|z|)^{\frac{1}{2}} \cdot (C_1)^{\frac{1}{2}} = \sqrt{C_1 C_2} |z|^{\frac{1}{2}}.$$
	Dividing both sides by $|z|^{\frac{1}{2}}$ directly yields \eqref{eq:3.3.4}$(iii)$.
	
	Furthermore, by the triangle inequality we obtain
	\begin{flalign}
		j^{1/2} |z|^{1/2} \left| \mathcal{A}_k(z)^{j-1} \mathcal{B}_k(z) \right| \le 
		j^{1/2} |z|^{1/2} \left(\left|e^{(j-1)z} \right| \left| \mathcal{B}_k(z) \right| + \left| \mathcal{A}_k(z)^{j-1} - e^{(j-1)z} \right|  \left| \mathcal{B}_k(z) \right|\right). \label{eq:3.3.7}
	\end{flalign}
	For $j\geq2, k\geq1, z\in [-L_k,0]$, we have
	\begin{flalign}
		j^{1/2} \left| z \right|^{1/2} |e^{(j-1)z} \mathcal{B}_k(z)|   &\le C j^{1/2} |z|^{1/2} e^{(j-1)z}\nonumber\\ 
		& = C (\frac{j}{j-1})^{1/2} (j-1)^{1/2} |z|^{1/2} e^{(j-1)z} \nonumber\\
		&\leq C  (\frac{j}{j-1})^{1/2} \leq C.\label{eq:3.3.8}
	\end{flalign}
	By \eqref{eq:3.3.1.1}, the case $j = 1 $ of \eqref{eq:3.3.8}  can be proved. Next, for $j\geq2, k\geq1, z\in [-L_k,0]$, combining with \eqref{eq:3.3.4}$(iv)$ and \eqref{eq:3.3.1.1} we obtain
	\begin{flalign}	
		j^{1/2} \left| z \right|^{1/2} \left| \mathcal{A}_k(z)^{j-1} - e^{(j-1)z} \right|  \left| \mathcal{B}_k(z) \right|   
		\le \frac{C}{j-1} j^{1/2} 
		\le \frac{C}{\sqrt{j-1}}  
		\le C  \label{eq:3.3.9}
	\end{flalign}
	The above estimate obviously holds for $j = 1$, and together with \eqref{eq:3.3.7}, \eqref{eq:3.3.8}, \eqref{eq:3.3.9} we obtain the proof of \eqref{eq:3.3.4}$(v)$.
\end{proof}

Finally, we state the main results of this paper.
\begin{theorem}[Main theorem]\label{main_th}
	Let $r\in[0,1)$ and $p\in[2,\infty)$.
	Suppose that Assumptions~\ref{non},
	\ref{ass:noise}, \ref{ass:init},
	and~\ref{ass:AB} hold.
	Let $k\geq1$ satisfy
	$\tau\lambda_{h,N_h}\leq L_k$.
	Then the approximation defined by
	\eqref{eq:3.2.1} satisfies
	\begin{equation}\label{eq:3.3.10}
		\|u(t_m)-u_{h,m}\|_{L^p(\Omega;H)}
		\leq C\bigl(h^{1+r}+\tau^{1/2}\bigr),
		\qquad m=0,1,\ldots,M.
	\end{equation}
	The constant $C$ is independent of $h$, $\tau$, $m$,
	and the number of internal stages $k$ satisfying
	the stability condition.
\end{theorem} 

\subsection{Error estimates for deterministic problems}
For the subsequent convergence proof, this section presents some error estimates for the fully discrete deterministic equation. For the sake of simplicity, we shall concentrate on the case \( A = -\Delta \) with \( \mathcal{D}(A) = H^2(D) \cap H_0^1(D) \), where \( D \subset \mathbb{R}^d \), \( d = 1, 2, 3 \) is a bounded convex polygonal domain.  First, we introduce such a deterministic problem
\begin{flalign}
	{du(t)} + Au(t) = F(u(t)), \quad t\in (0,T], \quad  u(0) = u_0, \label{eq:3.4.1}
\end{flalign}
the source term 
$F$ satisfies Assumption \ref{non}. Therefore, the mild solution of \eqref{eq:3.4.1} can be expressed as follows
\begin{flalign}
	u(t) = E(t)u_0 + \int_0^t E(t-s)F(u(s)) \, ds,  \quad t \in [0, T]. \label{eq:3.4.2}
\end{flalign}

To numerically solve the deterministic problem \eqref{eq:3.4.1}, we employ the same spatial finite element approximation and time discretization strategy as for the original stochastic problem. Utilizing the explicit stabilized operators $\mathcal{A}_k$ and $\mathcal{B}_k$ introduced in Section \ref{sec3.2}, we define the fully discrete approximation 
\begin{flalign}
	\begin{cases}
		u_{m,h} &= \mathcal{A}_k(-\tau A_h)^m u_{0,h} + \tau  \sum_{i=0}^{m-1} \mathcal{A}_k(-\tau A_h)^{m-1-i} \mathcal{B}_k(-\tau A_h) P_h F(u_{i,h}) \\[3pt]
		u_{0,h} &= P_hu_0,\quad m = 1,2,\cdots,M
		\label{eq:3.4.3}
	\end{cases}
\end{flalign}

The error analysis of the fully discrete scheme \eqref{eq:3.4.3} is divided into two parts: the approximation error associated with the homogeneous evolution operator and the error arising from the discrete convolution term. We first establish estimates for the homogeneous part. Similar deterministic error estimates are well known for the linear implicit Euler method; see, \cite{Kruse2014a,Wang2017}. In the present explicit stabilized setting, however, the analysis relies on the stability and regularization properties of $\mathcal{A}_k(-\tau A_h)$ and $\mathcal{B}_k(-\tau A_h)$established in Assumption~\ref{ass:AB} and Lemma~\ref{lemma3.1}. Compared with the linear implicit Euler method, the discrete evolution operator associated with the explicit stabilized method possesses weaker regularization properties. Therefore, the corresponding estimates cannot be obtained as a direct consequence of the classical implicit Euler analysis.

\begin{lemma}
	\label{th3}
	Let Assumption \ref{ass:AB} be satisfied, and let \(k\ge1\) be chosen such that
	$
	\tau\lambda_{h,N_h}\le L_k.
	$ Let \(t_m=m\tau\le T\). Then the following estimates hold.
	\begin{enumerate}
		\item For every \(\beta\in[0,1]\) and \(v\in\dot H^{2\beta}\),
		\begin{flalign}
			\|(E(t_m) - \mathcal{A}_k(-\tau A_h)^m P_h)v\| \leq C \bigl(h^{2\beta} + \tau^\beta \bigr) \|v\|_{2\beta}. \label{eq:3.4.4}
		\end{flalign}
		
		\item  For every \(v\in H\) and \(m=1,2,\ldots,M\),
		\begin{flalign}
			\|(E(t_m) - \mathcal{A}_k(-\tau A_h)^{m}P_h)v\| \leq C \bigl(h^2 + \tau\bigr) t_m^{-1} \|v\|. \label{eq:3.4.5}
		\end{flalign}
		Here \(C\) is independent of \(h,\tau,k\), and \(m\).	These estimates are crucial for the proof of the  Theorem \ref{main_th}.
	\end{enumerate}
\end{lemma}
\begin{proof}
	Using the triangle inequality
	\begin{flalign}
		\|(E(t_m) - \mathcal{A}_k(-\tau A_h)^m P_h)v\| 
		&\leq \|(E(t_m) - E_h(t_m) P_h)v\| + \|(E_h(t_m) - \mathcal{A}_k(-\tau A_h)^m) P_h v\| \nonumber \\
		&:= I_0 + I_1. \label{eq:3.4.6}
	\end{flalign}
	
	For \(\beta \in [0, 1]\) and \(v \in \dot{H}^{2\beta}\), by  \cite[Lemma~4.1(ii)]{Kruse2014a} we obtain
	\begin{flalign}
		I_0 \leq C h^{2\beta} \|v\|_{2\beta}.	\label{eq:3.4.7}
	\end{flalign}
	
	We denote by \( (\lambda_{h,n})_{n=1}^{N_h} \) the positive eigenvalues of \( A_h \) with corresponding orthonormal eigenvectors \( (\varphi_{h,n})_{n=1}^{N_h} \subset V_h \). For \(I_1\), using the spectral decomposition of \(A_h\), we obtain
	\begin{flalign}
		I_1^2 = \sum_{n=1}^{N_h} \lambda_{h,n}^{-2\beta} \left| e^{-m\tau \lambda_{h,n}} - \mathcal{A}_k(-\tau \lambda_{h,n})^m\right|^2	\left|(A_h^{\beta}P_h v, \varphi_{h,n})\right|^2.      \label{eq:3.4.8}  
	\end{flalign}  
	
	Therefore, for $k\geq 1$, we split the spectrum into the two index sets
	$$
	\mathcal{I}_{<}
	:=
	\left\{
	n\in\{1,\ldots,N_h\}:
	\tau\lambda_{h,n}\leq\gamma
	\right\},
	$$
	and
	$$
	\mathcal{I}_{>}
	:=
	\left\{
	n\in\{1,\ldots,N_h\}:
	\gamma<\tau\lambda_{h,n}\leq L_k
	\right\}.
	$$
	Hence,
	$$
	\{1,\ldots,N_h\}
	=
	\mathcal{I}_{<}\,\cup\,\mathcal{I}_{>}.
	$$
	
	For $\beta =0$, $m\geq1$,  combining with \eqref{eq:3.3.4}$(iv)$ and \eqref{eq:3.4.8}, we obtain
	\begin{flalign}
		I_1^2 &= \sum_{n=1}^{N_h} |e^{-m\tau \lambda_{h,n}} - \mathcal{A}_k(-\tau \lambda_{h,n})^m|^2 |(P_h v, \varphi_{n,h})|^2\nonumber\\
		&\leq C^2 m^{-2} \sum_{n=1}^{N_h} |(P_h v, \varphi_{n,h})|^2 \nonumber\\
		&\leq C^2 \sum_{n=1}^{N_h} |(P_h v, \varphi_{n,h})|^2 \nonumber\\
		&\leq C^2 \|v\|^2  
		\label{eq:3.4.9} 
	\end{flalign}
	
	When $\beta=1$, we get
	\begin{flalign}
		I_1^2 = \sum_{n=1}^{N_h} \lambda_{h,n}^{-2} |e^{-m\tau \lambda_{h,n}} - \mathcal{A}_k(-\tau \lambda_{h,n})^m |^2 (A_hP_h v, \varphi_{h,n})^2.  \label{eq:3.4.10} 
	\end{flalign}
	Next, define
	$$
	r_{1,k}(z)
	=
	\frac{e^{z}-\mathcal{A}_k(z)}{z^2}.
	$$
	By Assumption~\ref{ass:AB},
	$$
	\mathcal{A}_k(0)=1,
	\qquad
	\mathcal{A}_k'(0)=1.
	$$
	Thus the numerator $e^z-\mathcal{A}_k(z)$ has a zero of order at least two at
	$z=0$. Hence the singularity at $z=0$ is removable, and $r_{1,k}$ extends
	analytically to the open disc $\mathbb{D}(0,\frac{\delta}{2})$, where $\delta$ is given by
	\eqref{eq:3.3.3} and is independent of $k$.Owing to the maximum principle for analytic functions one has 
	\begin{equation}\label{eq:3.4.11}
		C_0
		:=
		\sup_{k\geq1,\ |z|\leq\delta/2}
		|r_{1,k}(z)|
		\leq
		\frac{4}{\delta^2}
		\sup_{k\geq1,\ |z|=\delta/2}
		\bigl(e^{\delta/2}+|\mathcal{A}_k(z)|\bigr)
		<\infty.
	\end{equation}
	Let $ \gamma \in (0, \delta) $ be chosen sufficiently small such that $ \frac{\gamma}{8} + C_0 \gamma e^{\gamma/2} < \frac{1}{2} $. Then, for all $ z \in [-\gamma, 0] $,  
	\begin{flalign}
		|\mathcal{A}_k(z)| &= e^{z/2} \left |e^{z/2} - e^{-z/2} \left(e^{z}-\mathcal{A}_k(z)\right) \right| \nonumber\\
		&\leq e^{z/2} (1 - |z|/2 + |z|^2/8 + C_0 |z|^2 e^{-z/2}) \leq e^{z/2}. \label{eq:3.4.12} 
	\end{flalign}
	For $n\in\mathcal{I}_{<}$, using
	$$
	e^{mz}-\mathcal{A}_k(z)^m
	=
	\bigl(e^z-\mathcal{A}_k(z)\bigr)
	\sum_{i=0}^{m-1}
	e^{iz}\mathcal{A}_k(z)^{m-1-i},
	$$
	together with
	$$
	e^z-\mathcal{A}_k(z)=z^2r_{1,k}(z)
	$$
	and \eqref{eq:3.4.12}, we obtain
	\begin{flalign}
		I_{1,<}^2
		&:=
		\sum_{n\in\mathcal{I}_{<}}
		\lambda_{h,n}^{-2}
		\left|
		e^{-m\tau\lambda_{h,n}}
		-
		\mathcal{A}_k(-\tau\lambda_{h,n})^m
		\right|^2
		\left|
		(A_hP_hv,\varphi_{h,n})
		\right|^2
		\nonumber\\
		&\leq
		\sum_{n\in\mathcal{I}_{<}}
		\lambda_{h,n}^{-2}
		|r_{1,k}(-\tau\lambda_{h,n})|^2
		(\tau\lambda_{h,n})^4
		\left|
		\sum_{i=0}^{m-1}
		e^{-i\tau\lambda_{h,n}}
		\mathcal{A}_k(-\tau\lambda_{h,n})^{m-1-i}
		\right|^2
		\left|
		(A_hP_hv,\varphi_{h,n})
		\right|^2
		\nonumber\\
		&\leq
		C_0^2\tau^4
		\sum_{n\in\mathcal{I}_{<}}
		\lambda_{h,n}^{2}
		m^2
		e^{-(m-1)\tau\lambda_{h,n}}
		\left|
		(A_hP_hv,\varphi_{h,n})
		\right|^2
		\nonumber\\
		&\leq
		C_0^2\tau^2 e^\gamma
		\sum_{n\in\mathcal{I}_{<}}
		(m\tau\lambda_{h,n})^2
		e^{-m\tau\lambda_{h,n}}
		\left|
		(A_hP_hv,\varphi_{h,n})
		\right|^2
		\nonumber\\
		&\leq
		C(\gamma)\tau^2
		\sum_{n\in\mathcal{I}_{<}}
		\left|
		(A_hP_hv,\varphi_{h,n})
		\right|^2,
		\label{eq:3.4.13}
	\end{flalign}
	where we have used $\sup_{y\geq0}y^2e^{-y}<\infty.$
	Next, for $n\in\mathcal{I}_{>}$, define
	$$
	\rho
	:=
	\sup_{k\geq1,z\in[-L_k,-\gamma]}
	|\mathcal{A}_k(z)|
	<1.
	$$
	Then
	\begin{flalign}
		I_{1,>}^2
		&:=
		\sum_{n\in\mathcal{I}_{>}}
		\lambda_{h,n}^{-2}
		\left|
		e^{-m\tau\lambda_{h,n}}
		-
		\mathcal{A}_k(-\tau\lambda_{h,n})^m
		\right|^2
		\left|
		(A_hP_hv,\varphi_{h,n})
		\right|^2
		\nonumber\\
		&\leq
		C
		\sum_{n\in\mathcal{I}_{>}}
		\lambda_{h,n}^{-2}
		\left(
		e^{-m\tau\lambda_{h,n}}
		+
		|\mathcal{A}_k(-\tau\lambda_{h,n})|^m
		\right)^2
		\left|
		(A_hP_hv,\varphi_{h,n})
		\right|^2
		\nonumber\\
		&\leq
		C
		\sum_{n\in\mathcal{I}_{>}}
		\lambda_{h,n}^{-2}
		\left(
		e^{-m\gamma}
		+
		\rho^m
		\right)^2
		\left|
		(A_hP_hv,\varphi_{h,n})
		\right|^2
		\nonumber\\
		&\leq
		C(\gamma)\tau^2
		\left(
		e^{-m\gamma}
		+
		\rho^m
		\right)^2
		\sum_{n\in\mathcal{I}_{>}}
		\left|
		(A_hP_hv,\varphi_{h,n})
		\right|^2
		\nonumber\\
		&\leq
		C(\gamma)\tau^2
		\sum_{n\in\mathcal{I}_{>}}
		\left|
		(A_hP_hv,\varphi_{h,n})
		\right|^2.
		\label{eq:3.4.14}
	\end{flalign}
	Here we have used the fact that
	$$
	\lambda_{h,n}^{-1}
	\leq
	\frac{\tau}{\gamma},
	\qquad
	n\in\mathcal{I}_{>}.
	$$
	Combining \eqref{eq:3.4.13} and \eqref{eq:3.4.14}, we obtain, for $\beta=1$,
	\begin{align}
		I_1^2
		&=
		I_{1,<}^2+I_{1,>}^2 \leq
		C\tau^2
		\sum_{n=1}^{N_h}
		\left|
		(A_hP_hv,\varphi_{h,n})
		\right|^2
		\nonumber\\
		&=
		C\tau^2\|A_hP_hv\|^2
		\leq
		C\tau^2\|v\|_2^2. \label{eq:I1}
	\end{align}
	Taking the square root gives $I_1 \le C \tau \|v\|_2$.
	Finally, interpolating between the results for $\beta=0$ in \eqref{eq:3.4.9}  and $\beta=1$, we deduce that for any $\beta \in [0, 1]$,
	\begin{flalign}I_1 \le C \tau^\beta \|v\|_{2\beta}. \label{eq:3.4.15}  
	\end{flalign}
	Substituting \eqref{eq:3.4.7} and \eqref{eq:3.4.15} into \eqref{eq:3.4.6}, we obtain \eqref{eq:3.4.4}. 
	
	Next, we prove the case of non-smooth initial data, and similarly we obtain
	\begin{flalign}
		\|(E(t_m) - \mathcal{A}_k(-\tau A_h)^m P_h)v\| 
		&\leq \|(E(t_m) - E_h(t_m) P_h)v\| + \|(E_h(t_m) - \mathcal{A}_k(-\tau A_h)^m) P_h v\| \nonumber \\
		&:= I_2 + I_3. \label{eq:3.4.16}
	\end{flalign}
	
	For \(v\in H\), the standard nonsmooth-data finite element estimate gives
	\begin{flalign}
		I_2 \le C h^2 t_m^{-1} \|v\|.	
		\label{eq:3.4.17}
	\end{flalign}
	
	For $ k\geq 1$, $m\geq1$, $\tau \lambda_{h,N_h} \in [0,L_k]$, combining with \eqref{eq:3.3.4}$(iv)$ and \eqref{eq:3.4.8}, we obtain
	\begin{flalign}
		I_3^2 &= \sum_{n=1}^{N_h} |e^{-m\tau \lambda_{h,n}} - \mathcal{A}_k(-\tau \lambda_{h,n})^m|^2 |(P_h v, \varphi_{h,n})|^2\nonumber\\
		&\leq C^2 m^{-2} \sum_{n=1}^{N_h} |(P_h v, \varphi_{h,n})|^2 \nonumber\\
		&= C^2 \tau^2 t_{m}^{-2} \sum_{n=1}^{N_h} |(P_h v, \varphi_{h,n})|^2 \nonumber\\
		&\leq C^2 \tau^2 t_{m}^{-2} \|v\|^2  
		\label{eq:3.4.18} 
	\end{flalign}
	Substituting the square roots of \eqref{eq:3.4.18} and \eqref{eq:3.4.17} into \eqref{eq:3.4.16} yields the proof for the non-smooth case. This completes the proof of Theorem~\ref{th3}.
\end{proof}
To estimate the error arising from the inhomogeneous term in (3.27), we next consider the discrete convolution operator. Setting \(j=m-i\), the relevant discrete propagator is
$$
A_k(-\tau A_h)^{j-1}B_k(-\tau A_h),\qquad j\ge1.
$$ The main task is to quantify its discrepancy from the continuous analytic semigroup \(E(t)\). Although the explicit stabilized evolution operator generally possesses weaker regularization properties than the linear implicit Euler propagator, the combined convolution operator still satisfies smoothing  of the same order needed in the convergence analysis. The corresponding estimates for \(v\in\dot H^{\mu-1}\) are stated in the following lemma and will subsequently be applied to the  source term \(F\).
\begin{lemma}
	\label{lemma3.6}
	Let Assumption \ref{ass:AB} be fulfilled and let \(k\ge1\) be chosen such that
	$
	\tau\lambda_{h,N_h}\le L_k$. For any $j = 1, 2, \ldots, M$, then it holds that
	\begin{flalign}
		\left\| \left( \mathcal{A}_k(-\tau A_h)^{j-1} \mathcal{B}_k(-\tau A_h) P_h - E(t) \right) v \right\| \leq C \left( h^\mu + \tau^{\frac{\mu}{2}} \right) t^{-\frac{1}{2}} \|v\|_{\mu-1}, \quad t \in (t_{j-1}, t_j],      	\label{eq:3.4.19} 
	\end{flalign}
	for all \( v \in \dot{H}^{\mu-1}, \, \mu \in [0, 2], \, h, \tau \in (0, 1] \).
\end{lemma}
\begin{proof}
	First consider \(\mu = 0\). Then we obtain the following estimate
	\begin{flalign}
		\|(A_k(-\tau A_h)^{j-1}B_k(-\tau A_h)P_h - E(t))v\| &\leq \|((A_k(-\tau A_h)^{j-1}B_k(-\tau A_h) - E_h(t_j))P_hv\| \nonumber\\
		&+ \|(E_h(t_j)P_h - E(t_j))v\| + \|(E(t_j) - E(t))v\| \nonumber\\
		&=: I_4 + I_5 + I_6. \label{eq:3.4.20}
	\end{flalign}
	
	Using Proposition \ref{pro1}, we obtain the following estimate
	\begin{flalign}	
		I_6 = \|(I - E(t_j - t))E(t)v\|
		&= \|(I - E(t_j - t)) A^{1/2}E(t) A^{-1/2}v\| \nonumber \\
		&\leq C  t^{-1/2}\|v\|_{-1}.	\label{eq:3.4.21}
	\end{flalign}
	
	For \(I_5\), using Lemma 4.1(ii) from \cite{Kruse2014a} with \(\rho = 1\), we obtain
	\begin{flalign}	
		I_5 \le  C t^{-1/2} \| v \|_{-1}. \label{eq:3.4.22}
	\end{flalign}
	
	For $I_4$, we can similarly obtain the following identity
	\begin{flalign}
		I_4^2 &= \|(\mathcal{A}_k(-\tau A_h)^{j-1}\mathcal{B}_k(-\tau A_h) - E_h(t_j))P_h v\|^2 \nonumber\\
		& = \sum_{n=1}^{N_h} \lambda_{h,n} \left| \mathcal{A}_k(-\tau \lambda_{h,n})^{j-1}\mathcal{B}_k(-\tau \lambda_{h,n}) - e^{-j\tau \lambda_{h,n}} \right|^2 (A_h^{-\frac{1}{2}} P_h v, \varphi_{h,n})^2. \label{eq:3.4.23}
	\end{flalign}
	For any $n\in\mathcal{I}_{<}$, $j\geq2$, by Lemma \ref{lemma3.1}, \eqref{eq:3.4.11} and \eqref{eq:3.4.12}, we have
	\begin{align}
		I_{4,<}^2
		&\le C \sum_{n\in\mathcal I_<} \lambda_{h,n}
		\left(
		\left|[\mathcal{A}_k(-\tau\lambda_{h,n})]^{j-1}
		\mathcal{B}_k(-\tau\lambda_{h,n})
		- e^{-(j-1)\tau\lambda_{h,n}}\mathcal{B}_k(-\tau\lambda_{h,n})\right|^2
		\right. \nonumber\\
		&\qquad\qquad\left.
		+ \left|e^{-(j-1)\tau\lambda_{h,n}}\mathcal{B}_k(-\tau\lambda_{h,n})
		- e^{-j\tau\lambda_{h,n}}\right|^2
		\right)
		(A_h^{-1/2}P_h v,\varphi_{h,n})^2 \nonumber\\
		&\le C \sum_{n\in\mathcal I_<} \lambda_{h,n}
		\left(
		C_0^2 j^2 e^{-(j-2)\tau\lambda_{h,n}}(\tau\lambda_{h,n})^4
		+ C_1^2(\tau\lambda_{h,n})^2 e^{-2(j-1)\tau\lambda_{h,n}}
		\right)
		(A_h^{-1/2}P_h v,\varphi_{h,n})^2 \nonumber\\
		&\le C(\gamma)(j\tau)^{-1}
		\sum_{n\in\mathcal I_<} j^{-2}
		\left(
		e^{-j\tau\lambda_{h,n}}(j\tau\lambda_{h,n})^5
		+ (j\tau\lambda_{h,n})^3 e^{-2j\tau\lambda_{h,n}}
		\right)
		(A_h^{-1/2}P_h v,\varphi_{h,n})^2 \nonumber\\
		&\le C(\gamma)(j\tau)^{-1}
		\sum_{n\in\mathcal I_<}
		(A_h^{-1/2}P_h v,\varphi_{h,n})^2. \label{eq:3.4.24}
	\end{align}
	For any $n\in\mathcal{I}_{>}$, $j\geq2$, by Lemma \ref{lemma3.1}, \eqref{eq:3.4.11} and \eqref{eq:3.4.12}, we have
	\begin{flalign}
		I_{4,>}^2 &\leq C	\sum_{n\in\mathcal{I}_{>}} \lambda_{h,n} 
		\Biggl( \Bigl| \mathcal{A}_k(-\tau \lambda_{h,n})^{j-1}\mathcal{B}_k(-\tau \lambda_{h,n})\Bigr|^2 +    
		\Bigl|e^{-j\tau \lambda_{h,n}} \Bigr|^2 \Biggl) (A_h^{-\frac{1}{2}} P_h v, \varphi_{h,n})^2 \nonumber \\
		&\leq C	\sum_{n\in\mathcal{I}_{>}} \lambda_{h,n}  	\left((j\tau\lambda_{h,n})^{-1} +(j\tau\lambda_{h,n})^{-1}(j\tau\lambda_{h,n}e^{-2j\tau \lambda_{h,n}})\right)(A_h^{-\frac{1}{2}} P_h v, \varphi_{h,n})^2\nonumber\\
		&\leq C	\sum_{n\in\mathcal{I}_{>}}	(j\tau)^{-1} \left(1+\left(j\tau\lambda_{h,n}e^{-2j\tau \lambda_{h,n}}\right)\right)(A_h^{-\frac{1}{2}} P_h v, \varphi_{h,n})^2\nonumber\\
		&\leq  C (j\tau)^{-1}  	\sum_{n\in\mathcal{I}_{>}}(A_h^{-\frac{1}{2}} P_h v, \varphi_{h,n})^2. \label{eq:3.4.25}
	\end{flalign} 
	When  $j = 1$, the above expression \eqref{eq:3.4.23} reduces to
	\begin{flalign}
		I_4^2 &= \sum_{n=1}^{N_h} \lambda_{h,n} \left| \mathcal{B}_k(-\tau \lambda_{h,n}) - e^{-\tau \lambda_{h,n}} \right|^2 (A_h^{-\frac{1}{2}} P_h v, \varphi_{h,n})^2 \nonumber\\
		& \leq C \sum_{n=1}^{N_h} \tau^{-1}   (\tau \lambda_{h,n}) \Biggl(\Bigl|\mathcal{B}_k(-\tau \lambda_{h,n}) \Bigl|^2 + \Bigl| e^{-\tau \lambda_{h,n}} \Bigl|^2 \Biggl)  (A_h^{-\frac{1}{2}} P_h v, \varphi_{h,n})^2, \nonumber\\
		& \leq C \tau^{-1} \sum_{n=1}^{N_h} (A_h^{-\frac{1}{2}} P_h v, \varphi_{h,n})^2 .
		\label{eq:3.4.26}
	\end{flalign}
	Combining the estimates for $I_4$, $I_5$ and $I_6$, the case $\mu=0$ can be proved. 
	
	Next, we prove the estimates for \(I_4\), \(I_5\), and \(I_6\) in the case \(\mu = 2\). Using Proposition \ref{pro1}, we obtain the following estimate 
	\begin{flalign}
		I_6 &= \|A^{-1}(I - E(t_j - t))A^{1/2}E(t) A^{1/2} v\| \nonumber\\
		&\le C    (t_j - t) t^{-1/2} \| v \|_{1} \nonumber\\
		&\le C   \tau  t^{-1/2} \| v \|_{1}   	 \label{eq:3.4.27}
	\end{flalign}
	
	For $I_5$, using Lemma 4.1(i) from \cite{Kruse2014a} with \(\mu = 2\), \(\nu = 1\), we obtain
	\begin{flalign}
		I_5 \le  C h^{2} t^{-1/2} \| v \|_{1}.  \label{eq:3.4.28}
	\end{flalign}
	
	For any $n\in\mathcal{I}_{<}$, $j\geq2$,  we adopt the same treatment as in the case $\mu = 2$
	\begin{align}
		I_{4,<}^2
		&= \sum_{n\in\mathcal I_<}\lambda_{h,n}^{-1}
		\left|\mathcal{A}_k(-\tau\lambda_{h,n})^{j-1}\mathcal{B}_k(-\tau\lambda_{h,n})
		-e^{-j\tau\lambda_{h,n}}\right|^2
		(A_h^{1/2}P_hv,\varphi_{h,n})^2 \nonumber\\
		&\le C\sum_{n\in\mathcal I_<}\lambda_{h,n}^{-1}
		\left(
		C_0^2j^2e^{-(j-2)\tau\lambda_{h,n}}(\tau\lambda_{h,n})^4
		+C_1^2(\tau\lambda_{h,n})^2e^{-2(j-1)\tau\lambda_{h,n}}
		\right)
		(A_h^{1/2}P_hv,\varphi_{h,n})^2 \nonumber\\
		&\le C(\gamma)\tau j^{-1}\sum_{n\in\mathcal I_<}
		\left(
		e^{-j\tau\lambda_{h,n}}(j\tau\lambda_{h,n})^3
		+e^{-2j\tau\lambda_{h,n}}(j\tau\lambda_{h,n})
		\right)
		(A_h^{1/2}P_hv,\varphi_{h,n})^2 \nonumber\\
		&\le C(\gamma)\tau^2t^{-1}\sum_{n\in\mathcal I_<}
		(A_h^{1/2}P_hv,\varphi_{h,n})^2.
		\label{eq:3.4.29}
	\end{align}
	For any $n\in\mathcal{I}_{>}$, $j\geq2$,
	\begin{flalign}
		I_{4,>}^2
		&\le C\sum_{n\in\mathcal I_>}\lambda_{h,n}^{-1}
		\left(
		\left|\mathcal{A}_k(-\tau\lambda_{h,n})^{j-1}\mathcal{B}_k(-\tau\lambda_{h,n})\right|^2
		+e^{-2j\tau\lambda_{h,n}}
		\right)
		(A_h^{1/2}P_hv,\varphi_{h,n})^2 \nonumber\\
		&\le C\sum_{n\in\mathcal I_>}\lambda_{h,n}^{-1}
		\left((j\tau\lambda_{h,n})^{-1}+e^{-2j\tau\lambda_{h,n}}\right)
		(A_h^{1/2}P_hv,\varphi_{h,n})^2 \nonumber\\
		&\le C\sum_{n\in\mathcal I_>}
		\left(
		\frac{1}{j\tau}\frac{\tau^2}{\gamma^2}
		+\frac{1}{j\tau}\frac{\tau^2}{\gamma^2}
		j\tau\lambda_{h,n}e^{-2j\tau\lambda_{h,n}}
		\right)
		(A_h^{1/2}P_hv,\varphi_{h,n})^2 \nonumber\\
		&\le C(\gamma)\tau^2t^{-1}
		\sum_{n\in\mathcal I_>}
		(A_h^{1/2}P_hv,\varphi_{h,n})^2.
		\label{eq:3.4.30}
	\end{flalign}
	When  $j = 1$, by \eqref{eq:3.3.4}$(iii)$ , we have
	\begin{flalign}
		I_4^2 &= \sum_{n=1}^{N_h} \lambda_{h,n}^{-1} \left| \mathcal{B}_k(-\tau \lambda_{h,n}) - e^{-\tau \lambda_{h,n}} \right|^2 (A_h^{\frac{1}{2}} P_h v, \varphi_{h,n})^2 \nonumber\\
		& \leq C \tau \sum_{n=1}^{N_h} (\tau \lambda_{h,n})^{-1} \left| \mathcal{B}_k(-\tau \lambda_{h,n}) - e^{-\tau \lambda_{h,n}} \right|^2 (A_h^{\frac{1}{2}} P_h v, \varphi_{h,n})^2, \nonumber\\
		& \leq C \tau \sum_{n=1}^{N_h} (A_h^{\frac{1}{2}} P_h v, \varphi_{h,n})^2.
		\label{eq:3.4.31}
	\end{flalign}
	Combining the estimates for $I_4$, $I_5$ and $I_6$, the case $\mu=2$ can be proved.  Finally, using an interpolation argument, we complete the proof of Theorem \ref{lemma3.6}. 	The subsequent theorem also plays a crucial role in the convergence analysis that follows.
\end{proof}

We also require the following error estimate for data of negative spatial regularity.

\begin{lemma}\label{lemma3.6zx}
	Let Assumption~\ref{ass:AB} hold, and let
	$h,\tau\in(0,1]$ and $k\ge1$ satisfy
	$\tau\lambda_{h,N_h}\le L_k$.
	Then, for every $\theta\in[0,1]$,
	$v\in\dot{H}^{-\theta}$, and $j=1,\ldots,M$,
	\begin{equation}\label{eq:negative-pointwise-error}
		\left\|
		\Bigl(
		\mathcal{A}_k(-\tau A_h)^{j-1}
		\mathcal{B}_k(-\tau A_h)P_h-E(t)
		\Bigr)v
		\right\|
		\le
		C\left(h^{2-\theta}+\tau^{(2-\theta)/2}\right)
		t^{-1}\|v\|_{-\theta},
		\qquad t\in(t_{j-1},t_j].
	\end{equation}
	The constant $C$ is independent of
	$h$, $\tau$, $k$, $j$, and $t$.
\end{lemma}

\begin{proof}
	Similarly, we first consider $\theta \in [0,1] $, 
	by interpolation between
	$\|P_hv\|\leq\|v\|$ and \eqref{eq3.1.2}, we have
	\begin{flalign}
		\|A_h^{-\theta/2}P_hv\|
		\leq C\|v\|_{-\theta},
		\qquad v\in\dot{H}^{-\theta}. \label{eq:pf_th3.6.1}
	\end{flalign}
	For $t\in(t_{j-1},t_j]$, then we obtain the following estimate
	\begin{flalign}
		\|(A_k(-\tau A_h)^{j-1}B_k(-\tau A_h)P_h - E(t))v\| &\leq \|((A_k(-\tau A_h)^{j-1}B_k(-\tau A_h) - E_h(t_j))P_hv\| \nonumber\\
		&+ \|(E_h(t_j)P_h - E(t_j))v\| + \|(E(t_j) - E(t))v\| \nonumber\\
		&=: \widetilde{I}_4 + \widetilde{I}_5 + \widetilde{I}_6. \label{eq:pf_th3.6.2}
	\end{flalign}
	
	For $\widetilde{I}_5$,   \cite[Lemma~4.1(iii)]{Kruse2014a}, applied to the finite element
	discretization considered here, gives
	$$
	\widetilde{I}_5
	\leq
	Ch^{2-\theta}t_j^{-1}\|v\|_{-\theta}
	\leq
	Ch^{2-\theta}t^{-1}\|v\|_{-\theta}. \label{eq:pf_th3.6.3}
	$$
	
	For $\widetilde{I}_6$, using Proposition \ref{pro1} and \eqref{eq:pf_th3.6.1}, we obtain the following estimate 
	\begin{flalign}
		\widetilde{I}_6
		&=
		\left\|
		A^{-(2-\theta)/2}
		(I-E(t_j-t))
		AE(t)A^{-\theta/2}v
		\right\|
		\nonumber \\
		&\leq
		C(t_j-t)^{(2-\theta)/2} t^{-1}
		\|v\|_{-\theta}
		\nonumber	\\
		&\leq
		C\tau^{(2-\theta)/2} t^{-1}
		\|v\|_{-\theta}. \label{eq:pf_th3.6.4}
	\end{flalign}
	
	Finally, we derive the estimate for $\widetilde{I}_4$
	in the case $\theta=0$.
	For $j\geq2$, using Lemma~\ref{lemma3.1}
	and \eqref{eq:3.3.1.1}, we obtain
	\begin{flalign}
		\widetilde{I}_4^2
		&=
		\sum_{n=1}^{N_h}
		\left|
		\mathcal{A}_k(-\tau\lambda_{h,n})^{j-1}
		\mathcal{B}_k(-\tau\lambda_{h,n})
		-
		e^{-j\tau\lambda_{h,n}}
		\right|^2
		|(P_hv,\varphi_{h,n})|^2
		\nonumber\\
		&\leq
		C\sum_{n=1}^{N_h}
		\Biggl(
		\left|
		\mathcal{A}_k(-\tau\lambda_{h,n})^{j-1}
		-
		e^{-(j-1)\tau\lambda_{h,n}}
		\right|^2
		|\mathcal{B}_k(-\tau\lambda_{h,n})|^2
		\nonumber\\
		&\qquad+
		e^{-2(j-1)\tau\lambda_{h,n}}
		\left|
		\mathcal{B}_k(-\tau\lambda_{h,n})
		-
		e^{-\tau\lambda_{h,n}}
		\right|^2
		\Biggr)
		|(P_hv,\varphi_{h,n})|^2
		\nonumber\\
		&\leq
		C\sum_{n=1}^{N_h}
		\left(
		(j-1)^{-2}
		+
		(\tau\lambda_{h,n})^2
		e^{-2(j-1)\tau\lambda_{h,n}}
		\right)
		|(P_hv,\varphi_{h,n})|^2
		\nonumber\\
		&=
		C(j-1)^{-2}
		\sum_{n=1}^{N_h}
		\left(
		1+
		((j-1)\tau\lambda_{h,n})^2
		e^{-2(j-1)\tau\lambda_{h,n}}
		\right)
		|(P_hv,\varphi_{h,n})|^2
		\nonumber\\
		&\leq
		Cj^{-2}\|P_hv\|^2
		\leq
		C\tau^2t_j^{-2}\|v\|^2
		\leq
		C\tau^2t^{-2}\|v\|^2,
		\label{eq:pf_th3.6.5}
	\end{flalign}
	here we used $\sup_{y\geq0}y^2e^{-2y}<\infty$.
	For $j=1$, \eqref{eq:3.3.1.1} and $t\leq\tau$
	give
	$\widetilde{I}_4^2\leq C\|v\|^2
	\leq C\tau^2t^{-2}\|v\|^2$.
	Therefore, $\widetilde{I}_4\leq C\tau t^{-1}\|v\|$
	for all $j\geq1$.
	
	We consider the case $\theta=1$.
	For $j\geq2$, using Lemma~\ref{lemma3.1},
	\eqref{eq:3.3.1.1} and \eqref{eq3.1.2}, we obtain
	\begin{flalign}
		\widetilde{I}_4^2
		&=
		\sum_{n=1}^{N_h}
		\lambda_{h,n}
		\left|
		\mathcal{A}_k(-\tau\lambda_{h,n})^{j-1}
		\mathcal{B}_k(-\tau\lambda_{h,n})
		-
		e^{-j\tau\lambda_{h,n}}
		\right|^2
		|(A_h^{-1/2}P_hv,\varphi_{h,n})|^2
		\nonumber\\
		&\leq
		C\sum_{n=1}^{N_h}
		\lambda_{h,n}
		\Biggl(
		\left|
		\mathcal{A}_k(-\tau\lambda_{h,n})^{j-1}
		-
		e^{-(j-1)\tau\lambda_{h,n}}
		\right|^2
		|\mathcal{B}_k(-\tau\lambda_{h,n})|^2
		\nonumber\\
		&\qquad+
		e^{-2(j-1)\tau\lambda_{h,n}}
		\left|
		\mathcal{B}_k(-\tau\lambda_{h,n})
		-
		e^{-\tau\lambda_{h,n}}
		\right|^2
		\Biggr)
		|(A_h^{-1/2}P_hv,\varphi_{h,n})|^2
		\nonumber\\
		&\leq
		C\tau^{-1}\sum_{n=1}^{N_h}
		\left(
		(j-1)^{-2}
		+
		(\tau\lambda_{h,n})^3
		e^{-2(j-1)\tau\lambda_{h,n}}
		\right)
		|(A_h^{-1/2}P_hv,\varphi_{h,n})|^2
		\nonumber\\
		&\leq
		C\tau^{-1}(j-1)^{-2}
		\sum_{n=1}^{N_h}
		\left(
		1+
		((j-1)\tau\lambda_{h,n})^3
		e^{-2(j-1)\tau\lambda_{h,n}}
		\right)
		|(A_h^{-1/2}P_hv,\varphi_{h,n})|^2
		\nonumber\\
		&\leq
		C\tau^{-1}j^{-2}\|A_h^{-1/2}P_hv\|^2
		\nonumber\\
		&\leq
		C\tau t_j^{-2}\|v\|_{-1}^2
		\nonumber\\
		&\leq
		C\tau t^{-2}\|v\|_{-1}^2,
		\label{eq:pf_th3.6.6}
	\end{flalign}
	here we used $\sup_{y\geq0}y^3e^{-2y}<\infty$.
	
	For $j=1$, \eqref{eq:3.3.1.1} gives
	\begin{flalign*}
		\widetilde{I}_4^2
		&=
		\sum_{n=1}^{N_h}
		\lambda_{h,n}
		\left|
		\mathcal{B}_k(-\tau\lambda_{h,n})
		-
		e^{-\tau\lambda_{h,n}}
		\right|^2
		|(A_h^{-1/2}P_hv,\varphi_{h,n})|^2
		\\
		&\leq
		C\tau^{-1}\sum_{n=1}^{N_h}
		\tau\lambda_{h,n}
		\left(
		|\mathcal{B}_k(-\tau\lambda_{h,n})|^2
		+
		e^{-2\tau\lambda_{h,n}}
		\right)
		|(A_h^{-1/2}P_hv,\varphi_{h,n})|^2
		\\
		&\leq
		C\tau^{-1}\|A_h^{-1/2}P_hv\|^2
		\\
		&\leq
		C\tau^{-1}\|v\|_{-1}^2
		\leq
		C\tau t^{-2}\|v\|_{-1}^2,
	\end{flalign*}
	Combining the estimates for $\widetilde{I}_4$, $\widetilde{I}_5$, and $\widetilde{I}_6$
	and interpolating between the cases $\theta=0$ and
	$\theta=1$ completes the proof.
\end{proof}

\begin{lemma}
	\label{lemma3.7}
	Let Assumption~\ref{ass:AB} be fulfilled and let
	\(\widetilde{\rho}\in[0,1]\).
	Let \(h,\tau\in(0,1]\), and let \(k\geq1\) be chosen such that
	$
	\tau\lambda_{h,N_h}\leq L_k.
	$
	Then, for every \(m=1,2,\ldots,M\) and every
	\(v\in\dot{H}^{-\widetilde{\rho}}\), there exists a constant
	\(C>0\), independent of \(h\), \(\tau\), \(k\), and \(m\), such that
	\begin{flalign}
		\left\|
		\sum_{j=1}^m
		\int_{t_{j-1}}^{t_j}
		\left(
		\mathcal{A}_k(-\tau A_h)^{j-1}
		\mathcal{B}_k(-\tau A_h)P_h
		-
		E(\sigma)
		\right)v\,d\sigma
		\right\|
		\leq
		C
		\left(
		h^{2-\widetilde{\rho}}
		+
		\tau^{(2-\widetilde{\rho})/2}
		\right)
		\|v\|_{-\widetilde{\rho}}.
		\label{eq:3.4.32}
	\end{flalign}
	Moreover, for any \( \beta \in [0, 2] \), \( v \in \dot{H}^{\beta-1} \) and \( h, \tau \in (0, 1] \), we have
	\begin{flalign}
		\left( \sum_{j=1}^m \int_{t_{j-1}}^{t_j}
		\left\| \left( [\mathcal{A}_k(-\tau A_h)]^{j-1}
		\mathcal{B}_k(-\tau A_h)P_h - E(\sigma) \right)v \right\|^2
		\mathrm{d}\sigma \right)^{1/2}
		\leq C \left( h^\beta + \tau^{\beta/2} \right) \|v\|_{\beta-1}.
		\label{eq:3.4.33}
	\end{flalign}
\end{lemma}
\begin{proof}
	We first prove estimate \eqref{eq:3.4.32}. Again we use the triangle inequality
	\begin{flalign}
		&\left\| \sum_{j=1}^m \int_{t_{j-1}}^{t_j}
		\left( [\mathcal{A}_k(-\tau A_h)]^{j-1}
		\mathcal{B}_k(-\tau A_h)P_h - E(\sigma) \right)v\,\mathrm{d}\sigma
		\right\| \nonumber\\
		&\leq \left\| \sum_{j=1}^m \int_{t_{j-1}}^{t_j}
		\left( [\mathcal{A}_k(-\tau A_h)]^{j-1}\mathcal{B}_k(-\tau A_h)
		- E_h(t_j) \right) P_h v \,\mathrm{d}\sigma \right\| \nonumber\\
		&\quad + \left\| \sum_{j=1}^m \int_{t_{j-1}}^{t_j}
		\left( E_h(t_j) - E_h(\sigma) \right) P_h v \,\mathrm{d}\sigma \right\| \nonumber\\
		&\quad + \left\| \sum_{j=1}^m \int_{t_{j-1}}^{t_j}
		\left( E_h(\sigma) P_h - E(\sigma) \right) v \,\mathrm{d}\sigma \right\| \nonumber\\
		&\quad := J_0 + J_1 + J_2 \label{eq:3.4.34}
	\end{flalign}
	
	According to \cite[Lemma~4.2(i)]{Kruse2014a}, we obtain the bound for 
	\begin{flalign}
		J_2 \leq Ch^{2-\widetilde{\rho}} \|v\|_{-\widetilde{\rho}}. \label{eq:3.4.35}
	\end{flalign}
	For convenience in the subsequent estimates, we introduce another partition of the index set, defined as follows
	$$
	\mathcal{J}_{<}
	:=
	\left\{
	n\in\{1,\ldots,N_h\}:
	\tau\lambda_{h,n}\leq 1
	\right\},
	$$
	and
	$$
	\mathcal{J}_{>}
	:=
	\left\{
	n\in\{1,\ldots,N_h\}:
	\tau\lambda_{h,n}>1
	\right\},
	$$
	where,
	$$
	\{1,\ldots,N_h\}
	=
	\mathcal{J}_{<}\,\cup\,\mathcal{J}_{>}.
	$$
	For any $n\in\mathcal{J}_{<}$, $j\geq1$, $\widetilde{\rho} \in \{0,1\}$,
	we use the  identity of $J_1$
	\begin{flalign}
		J_{1,<}^2
		&= \sum_{n\in\mathcal J_<}
		\left|
		\lambda_{h,n}^{\tilde\rho/2}
		\sum_{j=1}^m e^{-\lambda_{h,n}t_{j-1}}
		\int_0^\tau \left(e^{-\tau\lambda_{h,n}}-e^{-\lambda_{h,n}\sigma}\right)d\sigma
		\right|^2
		(A_h^{-\tilde\rho/2}P_hv,\varphi_{h,n})^2 \nonumber\\
		&= \sum_{n\in\mathcal J_<}
		\left|
		\lambda_{h,n}^{\tilde\rho/2}
		\left(\tau e^{-\tau\lambda_{h,n}}
		-\frac{1-e^{-\tau\lambda_{h,n}}}{\lambda_{h,n}}\right)
		\sum_{j=1}^m e^{-\tau\lambda_{h,n}(j-1)}
		\right|^2
		(A_h^{-\tilde\rho/2}P_hv,\varphi_{h,n})^2 \nonumber\\
		&\le C(\gamma)\sum_{n\in\mathcal J_<}
		\lambda_{h,n}^{\tilde\rho-2}
		e^{-2\tau\lambda_{h,n}}
		(\tau\lambda_{h,n})^4
		\left(1-e^{-\tau\lambda_{h,n}}\right)^{-2}
		(A_h^{-\tilde\rho/2}P_hv,\varphi_{h,n})^2 \nonumber\\
		&\le C(\gamma)\tau^{2-\tilde\rho}
		\sum_{n\in\mathcal J_<}
		(A_h^{-\tilde\rho/2}P_hv,\varphi_{h,n})^2 \nonumber\\
		&\le C(\gamma)\tau^{2-\tilde\rho}\|v\|_{-\tilde\rho}^2,
		\label{eq:3.4.36}
	\end{flalign}
	
	where we have used that the function \( x \mapsto x(1 - e^{-x})^{-1} \) is bounded for all \( x \in (0, 1] \). On the other hand, for  $n\in\mathcal{J}_{>}$, $j\geq1$, we have
	\begin{flalign}
		J_{1,>}^2 &= \sum_{n\in\mathcal{J}_{>}} \left| \lambda_{h,n}^{\frac{\widetilde{\rho}}{2}} \sum_{j=1}^m e^{-\lambda_{h,n} t_{j-1}} \int_0^\tau \left( e^{-\tau\lambda_{h,n}} - e^{-\lambda_{h,n} \sigma} \right) d\sigma \right|^2 \left( A_h^{-\frac{\widetilde{\rho}}{2}} P_h v, \varphi_{h,n} \right)^2 \nonumber \\ 
		&\leq  \sum_{n\in\mathcal{J}_{>}} 2 \lambda_{h,n}^{\widetilde{\rho}-2} \left( \left| \tau \lambda_{h,n} e^{-\tau\lambda_{h,n}} \right|^2 + \left| 1 - e^{-\tau\lambda_{h,n}} \right|^2 \right) \left( 1 - e^{-\tau\lambda_{h,n}} \right)^{-2}   \left( A_h^{-\frac{{\widetilde{\rho}}}{2}} P_h v, \varphi_{h,n} \right)^2 \nonumber \\ 
		&\leq C \tau^{2-\widetilde{\rho}} \sum_{n\in\mathcal{J}_{>}} \left( A_h^{-\frac{\widetilde{\rho}}{2}} P_h v, \varphi_{h,n} \right)^2 \nonumber \\ 
		&\leq C \tau^{2-\widetilde{\rho}} \left\|  v \right\|_{-\widetilde{\rho}}^2,  		\label{eq:3.4.37}
	\end{flalign}
	here we have used the $\lambda_{h,n}^{\widetilde{\rho}-2} < \tau^{2-\widetilde{\rho}} $ and $\left( 1 - e^{-\tau\lambda_{h,n}} \right)^{-2} \leq (1 - e^{-1})^{-2}$ inequalities.  
	For any $j\geq1$, $\widetilde{\rho} = 0$,$n\in\mathcal{J}_{<}$,
	\begin{flalign}
		J_0^2
		&= \sum_{n\in\mathcal J_<}
		\left|
		\tau\sum_{j=1}^m
		\left(\mathcal{A}_k(-\tau\lambda_{h,n})^{j-1}\mathcal{B}_k(-\tau\lambda_{h,n})
		-e^{-j\tau\lambda_{h,n}}\right)
		\right|^2
		(P_hv,\varphi_{h,n})^2 \nonumber\\
		&\le C(\gamma)\sum_{n\in\mathcal J_<}
		\left(
		\tau\sum_{j=1}^\infty
		\left(j e^{-j\tau\lambda_{h,n}/2}(\tau\lambda_{h,n})^2
		+(\tau\lambda_{h,n})e^{-j\tau\lambda_{h,n}}\right)
		\right)^2
		(P_hv,\varphi_{h,n})^2 \nonumber\\
		&\le C(\gamma)\tau^2\sum_{n\in\mathcal J_<}
		\left(
		(\tau\lambda_{h,n})^2
		\frac{e^{-\tau\lambda_{h,n}/2}}{(1-e^{-\tau\lambda_{h,n}/2})^2}
		+(\tau\lambda_{h,n})
		\frac{e^{-\tau\lambda_{h,n}}}{1-e^{-\tau\lambda_{h,n}}}
		\right)^2
		(P_hv,\varphi_{h,n})^2 \nonumber\\
		&\le C(\gamma)\tau^2\|v\|^2.
		\label{eq:3.4.38}
	\end{flalign}
	Since \( e^{-\tau\lambda_{h,n}} < 1 \), the conclusions concerning the geometric series in the above estimates are valid. Moreover, we also use the Taylor expansion \( 1 - e^{-x} = x + \mathcal{O}(x^2) \) as \( x \to 0 \). For any $j\geq1$, $\widetilde{\rho} = 0$, $n\in\mathcal{J}_{>}$, using Assumption \ref{ass:AB}
	\begin{flalign}
		J_0^2 &= \sum_{n\in\mathcal{J}_{>}}  \left| \tau \sum_{j=1}^m  \left( \mathcal{A}_k(-\tau\lambda_{h,n})^{j-1}\mathcal{B}_k(-\tau\lambda_{h,n}) - e^{-\lambda_{h,n}t_j} \right)  \right|^2 \left(  P_h v, \varphi_{h,n} \right)^2	\nonumber\\
		& \leq C \sum_{n\in\mathcal{J}_{>}}  \left( \tau \sum_{j=1}^m  \left| \mathcal{A}_k(-\tau\lambda_{h,n})^{j-1}\mathcal{B}_k(-\tau\lambda_{h,n}) - e^{-\lambda_{h,n}t_j} \right|   \right)^2 \left(  P_h v, \varphi_{h,n} \right)^2 \nonumber\\
		&\leq C \sum_{n\in\mathcal{J}_{>}}  \left( \tau \sum_{j=1}^m  \left| \mathcal{A}_k(-\tau\lambda_{h,n})^{j-1}\right| \left|\mathcal{B}_k(-\tau\lambda_{h,n})\right| + \left|e^{-\lambda_{h,n}t_j} \right|   \right)^2 \left(  P_h v, \varphi_{h,n} \right)^2 \nonumber\\
		&\leq C \sum_{n\in\mathcal{J}_{>}}  \left( \tau \sum_{j=1}^\infty \left( \rho^{j-1}  + e^{-\lambda_{h,n}t_j} \right) \right)^2 \left(  P_h v, \varphi_{h,n} \right)^2 \nonumber\\
		&\leq C \tau^2 \sum_{n\in\mathcal{J}_{>}} \left(  \frac{1}{1-\rho} + \frac{e^{-\tau\lambda_{h,n}}}{1-e^{-\tau\lambda_{h,n}}}  \right)^2 \left(  P_h v, \varphi_{h,n} \right)^2 \nonumber\\
		&\leq C(\gamma) \tau^2  \left\|  v \right\|^2.
		\label{eq:3.4.39}
	\end{flalign}
	
	Next we prove the case \(\widetilde{\rho} = 1\), for any $j\geq1$, $n\in\mathcal{I}_{<}$, by the same argument as for \(\widetilde{\rho} = 0\), we obtain
	\begin{flalign}
		J_{0,<}^2 &= \sum_{n\in\mathcal{I}_{<}}  \left| \lambda_{h,n}^{\frac{1}{2}} \sum_{j=1}^m \int_{t_{j-1}}^{t_j} \left( \mathcal{A}_k(-\tau\lambda_{h,n})^{j-1}\mathcal{B}_k(-\tau\lambda_{h,n}) - e^{-\lambda_{h,n}t_j} \right) d\sigma \right|^2 \left( A_h^{-\frac{1}{2}} P_h v, \varphi_{h,n} \right)^2	\nonumber\\
		&  \leq C(\gamma)\sum_{n\in\mathcal{I}_{<}} \left( \tau \lambda_{h,n}^{\frac{1}{2}} \sum_{j=1}^m  \left| \mathcal{A}_k(-\tau\lambda_{h,n})^{j-1}\mathcal{B}_k(-\tau\lambda_{h,n}) - e^{-\lambda_{h,n}t_j} \right|   \right)^2 \left( A_h^{-\frac{1}{2}} P_h v, \varphi_{h,n} \right)^2	\nonumber\\
		& \leq C(\gamma) \sum_{n\in\mathcal{I}_{<}} \left( \tau^{\frac{1}{2}} \sum_{j=1}^m \left(j e^{-\tau\lambda_{h,n}j/2} (\tau\lambda_{h,n})^2+(\tau\lambda_{h,n})e^{-j\tau\lambda_{h,n}}\right)   \right)^2  \left(A_h^{-\frac{1}{2}} P_h v, \varphi_{h,n} \right)^2	\nonumber\\
		&\leq C(\gamma) \tau \sum_{n\in\mathcal{I}_{<}}  \left((\tau\lambda_{h,n})^2 \frac{e^{-\tau\lambda_{h,n}/2}}{(1 - e^{-\tau\lambda_{h,n}/2})^2}+ (\tau\lambda_{h,n}) \frac{e^{-\tau\lambda_{h,n}}}{(1 - e^{-\tau\lambda_{h,n}})} \right)^2  \left(A_h^{-\frac{1}{2}} P_h v, \varphi_{h,n} \right)^2	 \nonumber\\
		&\leq C(\gamma) \tau \sum_{n\in\mathcal{I}_{<}} \left( A_h^{-\frac{1}{2}}P_h v, \varphi_{h,n} \right)^2	\nonumber\\
		&\leq C(\gamma) \tau  \left\|  v \right\|_{-1}^2.
		\label{eq:3.4.40}
	\end{flalign}
	For any  $j\geq1$, $n\in\mathcal{I}_{>}$,  using Assumption \ref{ass:AB} and Lemma \ref{lemma3.1}
	\begin{align}
		J_{0,>}^2
		&= \sum_{n\in\mathcal I_>}
		\left|
		\lambda_{h,n}^{1/2}\tau\sum_{j=1}^m
		\left(
		\mathcal{A}_k(-\tau\lambda_{h,n})^{j-1}\mathcal{B}_k(-\tau\lambda_{h,n})
		-e^{-j\tau\lambda_{h,n}}
		\right)
		\right|^2
		(A_h^{-1/2}P_hv,\varphi_{h,n})^2 \nonumber\\
		&\le C\sum_{n\in\mathcal I_>}
		\left(
		\tau^{1/2}\sum_{j=1}^m(\tau\lambda_{h,n})^{1/2}
		\left(
		\left|\mathcal{A}_k(-\tau\lambda_{h,n})^{j-1}\mathcal{B}_k(-\tau\lambda_{h,n})\right|
		+e^{-j\tau\lambda_{h,n}}
		\right)
		\right)^2
		(A_h^{-1/2}P_hv,\varphi_{h,n})^2 \nonumber\\
		&\le C\sum_{n\in\mathcal I_>}
		\left(
		\tau^{1/2}\sum_{j=1}^m
		\left(
		\rho^{j-1}
		+(\tau\lambda_{h,n})^{1/2}e^{-j\tau\lambda_{h,n}}
		\right)
		\right)^2
		(A_h^{-1/2}P_hv,\varphi_{h,n})^2 \nonumber\\
		&\le C\tau\sum_{n\in\mathcal I_>}
		\left(
		\frac{1}{1-\rho}
		+(\tau\lambda_{h,n})^{1/2}
		\frac{e^{-\tau\lambda_{h,n}}}{1-e^{-\tau\lambda_{h,n}}}
		\right)^2
		(A_h^{-1/2}P_hv,\varphi_{h,n})^2 \nonumber\\
		&\le C\tau\sum_{n\in\mathcal I_>}
		(A_h^{-1/2}P_hv,\varphi_{h,n})^2
		\le C\tau\|v\|_{-1}^2,
		\label{eq:3.4.41}
	\end{align}
	here we have used $\sup_{k \geq 1, z \in [-L_k, 0]}   |z|^{\frac{1}{2}}|\mathcal{B}_k(z)| \leq C$. In summary, we have completed the proof of \eqref{eq:3.4.32}. 
	
	For \eqref{eq:3.4.33}, we use the triangle inequality
	\begin{align}
		&\left( \sum_{j=1}^m \int_{t_{j-1}}^{t_j}
		\left\| \left( [\mathcal{A}_k(-\tau A_h)]^{j-1}
		\mathcal{B}_k(-\tau A_h)P_h - E(\sigma) \right)v \right\|^2
		\mathrm{d}\sigma \right)^{1/2}
		\nonumber\\
		&\le \left( \sum_{j=1}^m \int_{t_{j-1}}^{t_j}
		\left\| \left( [\mathcal{A}_k(-\tau A_h)]^{j-1}\mathcal{B}_k(-\tau A_h)
		- E_h(t_j) \right)P_hv \right\|^2
		\mathrm{d}\sigma \right)^{1/2}
		\nonumber\\
		&\quad + \left( \sum_{j=1}^m \int_{t_{j-1}}^{t_j}
		\left\| \left( E_h(t_j)-E_h(\sigma) \right)P_hv \right\|^2
		\mathrm{d}\sigma \right)^{1/2}
		\nonumber\\
		&\quad + \left( \sum_{j=1}^m \int_{t_{j-1}}^{t_j}
		\left\| \left( E_h(\sigma)P_h - E(\sigma) \right)v \right\|^2
		\mathrm{d}\sigma \right)^{1/2}
		\nonumber\\
		&\quad := J_3+J_4+J_5.
		\label{eq:3.4.42}
	\end{align}
	
	We first estimate $J_5$. For $\beta=0$, we initially take $v\in H$ and obtain
	\begin{flalign}
		J_5^2 &=  \sum_{j=1}^m \int_{t_{j-1}}^{t_j} \left\| \left( E_h(\sigma) P_h - E(\sigma) \right) v \right\|^{2} \, d\sigma 
		=  \int_0^{t_m} \left\| E_h(\sigma)P_h v - E(\sigma)v \right\|^2 d\sigma \ \nonumber \\
		&\le 2 \int_0^{t_m} \left\| E_h(\sigma)P_h v \right\|^2 d\sigma + 2 \int_0^{t_m} \| E(\sigma)v \|^2 d\sigma  \nonumber \\
		&\leq 2 \int_0^{t_m} \left( e^{-2\sigma A_h} P_h v, P_h v \right) d\sigma \ + 2\int_0^{t_m} \left( e^{-2\sigma A} v, v \right) d\sigma \ \nonumber \\
		&\leq\left( A_h^{-1}(I - e^{-2t_m A_h}) P_h v, P_h v \right) \ + \left( A^{-1}(I - e^{-2t_m A}) v, v \right) \   \nonumber \\
		&= \left\| A_h^{-\frac{1}{2}} P_h v \right\|^2 + \left\| A^{-\frac{1}{2}} v \right\|^2  \nonumber \\
		&\leq C \left\| v \right\|_{-1}^2.  \label{eq:3.4.43}
	\end{flalign}
	By density of $H$ in $\dot{H}^{-1}$, this estimate extends to every $v\in\dot{H}^{-1}$.
When $\beta = 2$, the same argument as in the proof of
\citet{Wang2017}, with the fully discrete propagator replaced by the
semigroup $E_h$, gives
$J_5 \leq C h^{2} \| v \|_{1}$.
	Therefore, for $\beta \in [0,2]$, combining  \eqref{eq:3.4.43}, we obtain the estimate of $J_5$
	\begin{flalign}
		J_5 \leq C h^{\beta} \| v \|_{\beta-1}.	\label{eq:3.4.44}
	\end{flalign}
	Next, we further estimate $J_4$. For $\beta = 2$, as in the proof of \eqref{eq:3.4.32} an application identity yields 
	\begin{flalign}
		J_4^2
		&= \sum_{j=1}^m\int_{t_{j-1}}^{t_j}
		\|(E_h(t_j)-E_h(\sigma))P_hv\|^2\mathrm{d}\sigma \nonumber\\
		&= \sum_{j=1}^m\sum_{n=1}^{N_h}\lambda_{h,n}^{-1}
		\int_{t_{j-1}}^{t_j}
		\left|e^{-j\tau\lambda_{h,n}}-e^{-\lambda_{h,n}\sigma}\right|^2\mathrm{d}\sigma
		(A_h^{1/2}P_hv,\varphi_{h,n})^2 \nonumber\\
		&\le \sum_{n=1}^{N_h}\lambda_{h,n}^{-1}\tau
		(1-e^{-\tau\lambda_{h,n}})^2
		(A_h^{1/2}P_hv,\varphi_{h,n})^2
		\sum_{j=1}^{\infty}e^{-2\tau\lambda_{h,n}(j-1)} \nonumber\\
		&= \sum_{n=1}^{N_h}\lambda_{h,n}^{-1}\tau
		\frac{1-e^{-\tau\lambda_{h,n}}}{1+e^{-\tau\lambda_{h,n}}}
		(A_h^{1/2}P_hv,\varphi_{h,n})^2 \nonumber\\
		&\le C\tau^2\|v\|_1^2.
		\label{eq:3.4.45}
	\end{flalign}
	
	On the other hand, for $\beta = 0$,  an identity yields  
	\begin{flalign}
		J_4^2&= \sum_{j=1}^{m} \sum_{n=1}^{N_h} \lambda_{h,n} \int_{t_{j-1}}^{t_{j}} \left| e^{-j \tau \lambda_{h,n}} - e^{- \lambda_{h,n} \sigma} \right|^2 \mathrm{d}\sigma \left( A_h^{-\frac{1}{2}} P_h v, \varphi_{h,n} \right)^2  \nonumber\\
		& =  \sum_{n=1}^{N_h} \lambda_{h,n} \left( \sum_{j=1}^m e^{-2(j-1)\tau\lambda_{h,n}} \right) \int_0^\tau \left( e^{-\tau\lambda_{h,n}} - e^{-\lambda_{h,n}s} \right)^2 ds \left( A_h^{-\frac{1}{2}} P_h v, \varphi_{h,n} \right)^2  \nonumber\\
		& =  \sum_{n=1}^{N_h} \lambda_{h,n} \left( \sum_{j=1}^m e^{-2(j-1)\tau\lambda_{h,n}} \right) 
		\left( \tau e^{-2\tau\lambda_{h,n}} - \frac{2e^{-\tau\lambda_{h,n}}}{\lambda_{h,n}}(1 - e^{-\tau\lambda_{h,n}}) + \frac{1 - e^{-2\tau\lambda_{h,n}}}{2\lambda_{h,n}} \right)   \left( A_h^{-\frac{1}{2}} P_h v, \varphi_{h,n}  \right)^2  \nonumber\\
		& \leq  \sum_{n=1}^{N_h} \lambda_{h,n}  \left(1- e^{-2\tau \lambda_{h,n}} \right)^{-1} 
		\left( \tau e^{-2\tau\lambda_{h,n}} - \frac{2e^{-\tau\lambda_{h,n}}}{\lambda_{h,n}}(1 - e^{-\tau\lambda_{h,n}}) + \frac{1 - e^{-2\tau\lambda_{h,n}}}{2\lambda_{h,n}} \right)   \left( A_h^{-\frac{1}{2}} P_h v, \varphi_{h,n}  \right)^2  \nonumber\\
		& \leq  \sum_{n=1}^{N_h} 
		\left(\frac{\tau\lambda_{h,n} e^{-2\tau\lambda_{h,n}}}{1 - e^{-2\tau\lambda_{h,n}}} + \frac{1 - 3e^{-\tau\lambda_{h,n}}}{2(1 + e^{-\tau\lambda_{h,n}})} \right)   \left( A_h^{-\frac{1}{2}} P_h v, \varphi_{h,n}  \right)^2  \nonumber\\
		& \leq C \sum_{n=1}^{N_h}  \left( A_h^{-\frac{1}{2}} P_h v, \varphi_{h,n}  \right)^2  \nonumber\\
		& \leq C \| v \|_{-1}^2,
		\label{eq:3.4.46}
	\end{flalign}
	here we use $ \sup_{x \geq 0} \left(f(x) = \frac{x e^{-2x}}{1 - e^{-2x}} + \frac{1 - 3e^{-x}}{2(1 + e^{-x})}\right) < \infty$.
	
	Finally, we aim to obtain an estimate for $J_3$, again giving priority to the case $\beta = 0$. For any  $j\geq1$, $n\in\mathcal{I}_{<}$,
	\begin{flalign}
		J_3^2&=  \sum_{n\in\mathcal{I}_{<}}  \lambda_{h,n} \sum_{j=1}^{m} \tau \left| \mathcal{A}_k(-\tau \lambda_{h,n})^{j-1} \mathcal{B}_k(-\tau\lambda_{h,n})-e^{-j \tau \lambda_{h,n}}   \right|^2 \left( A_h^{-\frac{1}{2}} P_h v, \varphi_{h,n} \right)^2  \nonumber\\
		&\leq C(\gamma)\sum_{n\in\mathcal{I}_{<}} \tau\lambda_{h,n} \sum_{j=1}^{m}  \left( j e^{-\tau\lambda_{h,n}j/2} (\tau\lambda_{h,n})^2+(\tau\lambda_{h,n})e^{-j\tau\lambda_{h,n}}  \right)^2 \left( A_h^{-\frac{1}{2}} P_h v, \varphi_{h,n} \right)^2  \nonumber\\
		&\leq C(\gamma) \sum_{n\in\mathcal{I}_{<}} \left(\sum_{j=1}^{m}  \left( j e^{-\tau\lambda_{h,n}j/2} (\tau\lambda_{h,n})^2+(\tau\lambda_{h,n})e^{-j\tau\lambda_{h,n}} \right) \right)^2 \left( A_h^{-\frac{1}{2}} P_h v, \varphi_{h,n} \right)^2  \nonumber\\
		&\leq C(\gamma)  \sum_{n\in\mathcal{I}_{<}}  \left((\tau\lambda_{h,n})^2 \frac{e^{-\tau\lambda_{h,n}/2}}{(1 - e^{-\tau\lambda_{h,n}/2})^2}+ (\tau\lambda_{h,n}) \frac{e^{-\tau\lambda_{h,n}}}{(1 - e^{-\tau\lambda_{h,n}})} \right)^2  \left( A_h^{-\frac{1}{2}} P_h v, \varphi_{h,n} \right)^2 \nonumber\\
		&\leq C(\gamma) \sum_{n\in\mathcal{I}_{<}}  \left( A_h^{-\frac{1}{2}}P_h v, \varphi_{h,n} \right)^2	\nonumber\\
		& \leq  C(\gamma) \| v \|_{-1}^2.	\label{eq:3.4.47}
	\end{flalign}
	For any  $j\geq1$, $n\in\mathcal{I}_{>}$, using  $\sup_{k \geq 1, z \in [-L_k, 0]}   |z|^{\frac{1}{2}}|\mathcal{B}_k(z)| \leq C$
	\begin{flalign}
		J_3^2
		&= \sum_{n\in\mathcal I_>}\lambda_{h,n}\tau\sum_{j=1}^m
		\left|
		\mathcal{A}_k(-\tau\lambda_{h,n})^{j-1}\mathcal{B}_k(-\tau\lambda_{h,n})
		-e^{-j\tau\lambda_{h,n}}
		\right|^2
		(A_h^{-1/2}P_hv,\varphi_{h,n})^2 \nonumber\\
		&\le C\sum_{n\in\mathcal I_>}\lambda_{h,n}\tau\sum_{j=1}^m
		\left(
		\left|\mathcal{A}_k(-\tau\lambda_{h,n})^{j-1}\mathcal{B}_k(-\tau\lambda_{h,n})\right|
		+e^{-j\tau\lambda_{h,n}}
		\right)^2
		(A_h^{-1/2}P_hv,\varphi_{h,n})^2 \nonumber\\
		&\le C\sum_{n\in\mathcal I_>}
		\left(
		\sum_{j=1}^m
		\left(
		\rho^{j-1}
		+(\tau\lambda_{h,n})^{1/2}e^{-j\tau\lambda_{h,n}}
		\right)
		\right)^2
		(A_h^{-1/2}P_hv,\varphi_{h,n})^2 \nonumber\\
		&\le C\sum_{n\in\mathcal I_>}
		\left(
		\frac{1}{1-\rho}
		+(\tau\lambda_{h,n})^{1/2}
		\frac{e^{-\tau\lambda_{h,n}}}{1-e^{-\tau\lambda_{h,n}}}
		\right)^2
		(A_h^{-1/2}P_hv,\varphi_{h,n})^2 \nonumber\\
		&\le C\sum_{n\in\mathcal I_>}
		(A_h^{-1/2}P_hv,\varphi_{h,n})^2
		\le C\|v\|_{-1}^2.
		\label{eq:3.4.48}
	\end{flalign}
	On the other hand, for $\beta = 2$, $j\geq1$,  $n\in\mathcal{I}_{<}$,   we can similarly obtain 
	\begin{flalign}
		J_3^2&=  \sum_{n\in\mathcal{I}_{<}} \lambda_{h,n}^{-1} \sum_{j=1}^{m} \tau \left| \mathcal{A}_k(-\tau \lambda_{h,n})^{j-1} \mathcal{B}_k(-\tau\lambda_{h,n})-e^{-j \tau \lambda_{h,n}}   \right|^2 \left( A_h^{\frac{1}{2}} P_h v, \varphi_{h,n} \right)^2  \nonumber\\
		&\le C(\gamma)  \sum_{n\in\mathcal{I}_{<}}  \lambda_{h,n}^{-1} \tau \sum_{j=1}^m \left( j e^{-\tau\lambda_{h,n}j/2} (\tau\lambda_{h,n})^2 + (\tau\lambda_{h,n}) e^{-j\tau\lambda_{h,n}} \right)^2 \left( A_h^{\frac{1}{2}} P_h v, \varphi_{h,n} \right)^2 \nonumber\\
		&\le C(\gamma)  \sum_{n\in\mathcal{I}_{<}}  \lambda_{h,n}^{-1} \tau \left( \sum_{j=1}^m j^2 e^{-j\tau\lambda_{h,n}} (\tau\lambda_{h,n})^4 + \sum_{j=1}^m e^{-2j\tau\lambda_{h,n}} (\tau\lambda_{h,n})^2 \right) \left( A_h^{\frac{1}{2}} P_h v, \varphi_{h,n} \right)^2 \nonumber\\
		&\le C(\gamma) \sum_{n\in\mathcal{I}_{<}}  \tau\lambda_{h,n}^{-1}  \left( (\tau\lambda_{h,n})^4  \frac{e^{-\tau\lambda_{h,n}}(1+e^{-\tau\lambda_{h,n}})}{(1-e^{-\tau\lambda_{h,n}})^3} + (\tau\lambda_{h,n})^2 \frac{e^{-2\tau\lambda_{h,n}}}{(1 - e^{-2\tau\lambda_{h,n}})} \right) \left( A_h^{\frac{1}{2}} P_h v, \varphi_{h,n} \right)^2 \nonumber\\
		&\le C(\gamma) \sum_{n\in\mathcal{I}_{<}}  \tau\lambda_{h,n}^{-1}  \left( (\tau\lambda_{h,n})^4 \frac{1}{ (\tau\lambda_{h,n})^3} + (\tau\lambda_{h,n})^2 \frac{1}{2 (\tau\lambda_{h,n})}  \right) \left( A_h^{\frac{1}{2}} P_h v, \varphi_{h,n} \right)^2 \nonumber\\
		&\leq \frac{3}{2} C(\gamma) \tau^{2}  \sum_{n\in\mathcal{I}_{<}}  \left( A_h^{\frac{1}{2}}P_h v, \varphi_{h,n} \right)^2	\nonumber\\
		& \leq  C(\gamma) \tau^{2} \| v \|_{1}^2,	\label{eq:3.4.49}
	\end{flalign}
	here Taylor expansion \( 1 - e^{-x} = x + \mathcal{O}(x^2) \) as \( x \to 0 \) is also used. For $\beta = 2$,  $j\geq1$,  $n\in\mathcal{I}_{>}$,  we can  obtain  
	\begin{flalign}
		J_3^2&=  \sum_{n\in\mathcal{I}_{>}} \lambda_{h,n}^{-1} \sum_{j=1}^{m} \tau \left| \mathcal{A}_k(-\tau \lambda_{h,n})^{j-1} \mathcal{B}_k(-\tau\lambda_{h,n})-e^{-j \tau \lambda_{h,n}}   \right|^2 \left( A_h^{\frac{1}{2}} P_h v, \varphi_{h,n} \right)^2  \nonumber\\
		&\leq C  \sum_{n\in\mathcal{I}_{>}} \lambda_{h,n}^{-1}\tau \sum_{j=1}^{m}  \left( \left| \mathcal{A}_k(-\tau\lambda_{h,n})^{j-1}\right|\left|\mathcal{B}_k(-\tau\lambda_{h,n})\right| + \left|e^{-\lambda_{h,n}t_j} \right| \right)^2 \left( A_h^{\frac{1}{2}} P_h v, \varphi_{h,n} \right)^2  \nonumber\\
		&\leq C  \sum_{n\in\mathcal{I}_{>}} \lambda_{h,n}^{-1}\tau \left(\sum_{j=1}^{\infty}  \left( \left| \mathcal{A}_k(-\tau\lambda_{h,n})^{j-1}\right|+ \left|e^{-j\gamma} \right| \right) \right)^2 \left( A_h^{\frac{1}{2}} P_h v, \varphi_{h,n} \right)^2  \nonumber\\
		&\leq C  \sum_{n\in\mathcal{I}_{>}} \lambda_{h,n}^{-1}\tau \left( \frac{1}{1-\rho} + \frac{e^{-\gamma}}{1-e^{-\gamma}} \right)^2 \left( A_h^{\frac{1}{2}} P_h v, \varphi_{h,n} \right)^2  \nonumber\\
		&\le C(\gamma) \sum_{n\in\mathcal{I}_{>}} (\lambda_{h,n}^{-1} \tau) \left( A_h^{\frac{1}{2}} P_h v, \varphi_{h,n} \right)^2\nonumber\\
		&\le C(\gamma) \frac{\tau^2}{\gamma} \sum_{n\in\mathcal{I}_{>}} \left( A_h^{\frac{1}{2}} P_h v, \varphi_{h,n} \right)^2 \nonumber\\
		&\le C(\gamma) \tau^2 \| v \|_1^2. \label{eq:3.4.50}
	\end{flalign}
	Combining the above estimates and using interpolation techniques, for $\beta \in [0,2]$, There exists a constant $c$ independent of $\tau,k$ such that
	\begin{flalign}
		J_3 \leq C \tau^{\frac{\beta}{2}} \| v \|_{\beta-1}.	\label{eq:3.4.51}
	\end{flalign}
	Combining \eqref{eq:3.4.44}, \eqref{eq:3.4.45}, \eqref{eq:3.4.46} and \eqref{eq:3.4.51}, we prove that \eqref{eq:3.4.33}.
	Therefore, we have completed the proof of the theorem \ref{lemma3.7}.
\end{proof}
In summary, although the discrete evolution operators associated with the explicit stabilized method possess weaker regularization properties than the linear implicit Euler propagator, the estimates derived above show that the corresponding fully discrete convolution operators achieve the same  smoothing required in the implicit Euler analysis. These estimates will be used in the subsequent strong convergence analysis.

\section{Proof of Theorem \ref{main_th}} \label{sec4}
This section is devoted to the proof of
Theorem~\ref{main_th}.
The proof relies on the Burkholder--Davis--Gundy
inequality and a discrete Gronwall inequality
with a weakly singular kernel.
We begin with a variant of the discrete Gronwall
lemma in \cite[Lemma~7.1]{ElliottLarsson1992}.

\begin{lemma}
	\label{lem:4.1}
	Let $T>0$, $\tau>0$, and $M\in\mathbb{N}$ satisfy
	$M\tau\le T$, and set $t_m=m\tau$.
	Let $C_1,C_2\ge0$, $\beta\in(0,1]$, and
	$(\varphi_m)_{m=1}^{M}$ be a nonnegative sequence
	satisfying
	\begin{equation}\label{eq:4.1}
		\varphi_m
		\le C_1
		+C_2\tau\sum_{i=1}^{m-1}
		(t_m-t_i)^{\beta-1}\varphi_i,
		\qquad m=1,2,\ldots,M.
	\end{equation}
	Then there exists a constant $C=C(C_2,T,\beta)$
	such that
	$$
	\varphi_m\le CC_1,
	\qquad m=1,2,\ldots,M.
	$$
	In particular, $C$ is independent of $\tau$ and $M$.
\end{lemma}

\begin{proof}
	See \cite[Lemma~6.1]{Kruse2014a}, which gives a variant
	of \cite[Lemma~7.1]{ElliottLarsson1992}.
\end{proof}
To estimate the stochastic integral terms, we recall the following Burkholder--Davis--Gundy inequality.
\begin{lemma}
	\label{lem:bdg}
	Let $p\in[2,\infty)$, and let
	$\Phi:[0,T]\times\Omega\to\mathcal{L}_2^0$
	be a predictable process satisfying
	$
	\mathbb{E}\left[
	\left(
	\int_0^T
	\|\Phi(s)\|_{\mathcal{L}_2^0}^2
	\,\mathrm{d}s
	\right)^{p/2}
	\right]<\infty.
	$
	Then, for every $t\in[0,T]$,
	\begin{equation}\label{eq:bdg}
		\left\|
		\int_0^t\Phi(s)\,\mathrm{d}W^Q(s)
		\right\|_{L^p(\Omega;H)}
		\leq
		C_p\left\|
		\left(
		\int_0^t
		\|\Phi(s)\|_{\mathcal{L}_2^0}^2
		\,\mathrm{d}s
		\right)^{1/2}
		\right\|_{L^p(\Omega;\mathbb{R})},
	\end{equation}
	where $C_p$ depends only on $p$.
\end{lemma}
\begin{proof}
	See \cite[Lemma~5.1]{Kruse2014a}.
\end{proof}
\begin{proof}[Proof of Theorem \ref{main_th}]
	For convenience, we write
	$$
	\mathcal{K}_j
	= \mathcal{A}_k(-\tau A_h)^{j-1} \mathcal{B}_k(-\tau A_h)P_h,
	\qquad j\geq1.
	$$
	Using the mild formulation \eqref{eq:1.2}
	and the discrete mild formulation of
	\eqref{eq:3.2.3}, with $u_{0,h}=P_hu_0$, we obtain
	\begin{align}
		&\|u(t_m)-u_{h,m}\|_{L^p(\Omega;H)}
		\le \left\| \bigl(E(t_m)-\mathcal{A}_k(-\tau A_h)^mP_h\bigr)u_0 \right\|_{L^p(\Omega;H)}	
		\notag\\
		&\quad+\Bigg\|	\sum_{i=0}^{m-1}	\int_{t_i}^{t_{i+1}}	\Bigl[	E(t_m-s)F(u(s))
		-\mathcal{K}_{m-i} F(u_{i,h})
		\Bigr]\,\mathrm{d}s	\Bigg\|_{L^p(\Omega;H)}
		\notag\\
		&\quad+\Bigg\|	\sum_{i=0}^{m-1}
		\int_{t_i}^{t_{i+1}}
		\Bigl[	E(t_m-s)G(u(s))	-\mathcal{K}_{m-i}G(u_{i,h})	\Bigr]\,\mathrm{d}W^Q(s)
		\Bigg\|_{L^p(\Omega;H)}
		\notag\\
		&:=\epsilon_0+\epsilon_1+\epsilon_2.
		\label{eq:4.3}
	\end{align}
	
	We estimate these three terms separately. For the initial-value error, applying \eqref{eq:3.4.4} with $\beta=(1+r)/2$ yields
	\begin{equation}\label{eq:4.4}
		\epsilon_0
		\le
		C\left(h^{1+r}+\tau^{(1+r)/2}\right)
		\|u_0\|_{L^p(\Omega;\dot{H}^{1+r})}.
	\end{equation}
	To estimate $\epsilon_1$, we use the following decomposition
	\begin{flalign}
		\epsilon_1
		&\le
		\Bigg\|
		\sum_{i=0}^{m-1}
		\int_{t_i}^{t_{i+1}}
		\mathcal{K}_{m-i}
		\bigl(F(u(s))-F(u_{i,h})\bigr)
		\,\mathrm{d}s
		\Bigg\|_{L^p(\Omega;H)}
		\notag\\
		&\quad+
		\Bigg\|
		\sum_{i=0}^{m-1}
		\int_{t_i}^{t_{i+1}}
		\Bigl(
		E(t_m-s)
		-\mathcal{K}_{m-i}
		\Bigr)
		\bigl(F(u(s))-F(u(t_m))\bigr)
		\,\mathrm{d}s
		\Bigg\|_{L^p(\Omega;H)}
		\notag\\
		&\quad+
		\Bigg\|
		\sum_{i=0}^{m-1}
		\int_{t_i}^{t_{i+1}}
		\Bigl(
		E(t_m-s)
		-\mathcal{K}_{m-i}
		\Bigr)
		F(u(t_m))
		\,\mathrm{d}s
		\Bigg\|_{L^p(\Omega;H)}
		\notag\\
		&:=\epsilon_{11}+\epsilon_{12}+\epsilon_{13}.
		\label{eq:4.5}
	\end{flalign}
	
	For $\epsilon_{11}$, interpolating the $H$-stability
	bound with the $\dot{H}^{-1}\to H$ estimate obtained
	from \eqref{eq:3.3.4}$(v)$ yields
	\begin{equation}\label{eq:drift-discrete-smoothing}
		\left\|
		\mathcal{K}_{m-i}v
		\right\|
		\le
		C t_{m-i}^{-(1-r)/2}\|v\|_{r-1},
		\qquad v\in\dot{H}^{r-1}.
	\end{equation}
	Using \eqref{eq:drift-discrete-smoothing},
	Assumption~\ref{non}, and the temporal regularity
	in Theorem~\ref{th1}, we obtain
	\begin{align}
		\epsilon_{11}
		&\le
		\sum_{i=0}^{m-1}
		\int_{t_i}^{t_{i+1}}
		\left\|
		\mathcal{K}_{m-i} P_h
		\bigl(F(u(s))-F(u_{i,h})\bigr)
		\right\|_{L^p(\Omega;H)}
		\,\mathrm{d}s
		\notag\\
		&\le
		C\sum_{i=0}^{m-1}
		t_{m-i}^{-(1-r)/2}
		\int_{t_i}^{t_{i+1}}
		\|F(u(s))-F(u_{i,h})\|_{L^p(\Omega;\dot{H}^{r-1})}
		\,\mathrm{d}s
		\notag\\
		&\le
		C\sum_{i=0}^{m-1}
		t_{m-i}^{-(1-r)/2}
		\int_{t_i}^{t_{i+1}}
		\|u(s)-u_{i,h}\|_{L^p(\Omega;H)}
		\,\mathrm{d}s
		\notag\\
		&\le
		C\tau\sum_{i=0}^{m-1}
		t_{m-i}^{-(1-r)/2}
		\|u(t_i)-u_{i,h}\|_{L^p(\Omega;H)}
		+
		C\sum_{i=0}^{m-1}
		t_{m-i}^{-(1-r)/2}
		\int_{t_i}^{t_{i+1}}(s-t_i)^{1/2}\,\mathrm{d}s
		\notag\\
		&\le
		C\tau\sum_{i=0}^{m-1}
		t_{m-i}^{-(1-r)/2}
		\|u(t_i)-u_{i,h}\|_{L^p(\Omega;H)}
		+
		C\tau^{1/2}
		\int_0^{t_m}s^{-(1-r)/2}\,\mathrm{d}s
		\notag\\
		&\le
		C\tau\sum_{i=0}^{m-1}
		t_{m-i}^{-(1-r)/2}
		\|u(t_i)-u_{i,h}\|_{L^p(\Omega;H)}
		+C\tau^{1/2}.
		\label{eq:4.6}
	\end{align}
	For $\epsilon_{12}$, we use the error estimate of Lemma~\ref{lemma3.6zx},
	together with Assumption~\ref{non}
	and Theorem~\ref{th1}, this yields
	\begin{align}
		\epsilon_{12}
		&\le
		\sum_{i=0}^{m-1}
		\int_{t_i}^{t_{i+1}}
		\Bigl\|
		\Bigl(
		E(t_m-s)
		-\mathcal{K}_{m-i}P_h
		\Bigr)
		\bigl(F(u(s))-F(u(t_m))\bigr)
		\Bigr\|_{L^p(\Omega;H)}
		\,\mathrm{d}s
		\notag\\
		&\le
		C\left(h^{1+r}+\tau^{(1+r)/2}\right)
		\int_0^{t_m}
		(t_m-s)^{-1}
		\|F(u(s))-F(u(t_m))\|_{L^p(\Omega;\dot{H}^{r-1})}
		\,\mathrm{d}s
		\notag\\
		&\le
		C\left(h^{1+r}+\tau^{(1+r)/2}\right)
		\int_0^{t_m}
		(t_m-s)^{-1}
		\|u(s)-u(t_m)\|_{L^p(\Omega;H)}
		\,\mathrm{d}s
		\notag\\
		&\le
		C\left(h^{1+r}+\tau^{(1+r)/2}\right)
		\int_0^{t_m}(t_m-s)^{-1/2}\,\mathrm{d}s
		\notag\\
		&\le
		C\left(h^{1+r}+\tau^{(1+r)/2}\right).
		\label{eq:4.7}
	\end{align}
	
	For $\epsilon_{13}$, applying \eqref{eq:3.4.32}
	with $\widetilde{\rho}=1-r$, together with
	Assumption~\ref{non} and Theorem~\ref{th1}, we obtain
	\begin{flalign}
		\epsilon_{13}
		&\leq
		C\left(h^{1+r}+\tau^{(1+r)/2}\right)
		\|F(u(t_m))\|_{L^p(\Omega;\dot{H}^{r-1})}
		\nonumber\\
		&\leq
		C\left(h^{1+r}+\tau^{(1+r)/2}\right)
		\left(
		1+
		\sup_{0\leq s\leq T}
		\|u(s)\|_{L^p(\Omega;H)}
		\right)
		\nonumber\\
		&\leq
		C\left(h^{1+r}+\tau^{(1+r)/2}\right).
		\label{eq:4.8}
	\end{flalign}

	Finally, we estimate the stochastic error term $\epsilon_2$, by the Burkholder--Davis--Gundy inequality and the triangle inequality, we obtain
	\begin{flalign}
		\epsilon_2
		&\leq
		C_p\left\|
		\left(
		\sum_{i=0}^{m-1}
		\int_{t_i}^{t_{i+1}}
		\left\|
		E(t_m-s)G(u(s))
		-
		\mathcal{K}_{m-i}G(u_{i,h})
		\right\|_{\mathcal{L}_2^0}^{2}
		\,\mathrm{d}s
		\right)^{1/2}
		\right\|_{L^p(\Omega;\mathbb{R})}
		\nonumber\\
		&\leq
		C_p\left\|
		\left(
		\sum_{i=0}^{m-1}
		\int_{t_i}^{t_{i+1}}
		\left\|
		\mathcal{K}_{m-i}
		\bigl(G(u(s))-G(u_{i,h})\bigr)
		\right\|_{\mathcal{L}_2^0}^{2}
		\,\mathrm{d}s
		\right)^{1/2}
		\right\|_{L^p(\Omega;\mathbb{R})}
		\nonumber\\
		&+
		C_p\left\|
		\left(
		\sum_{i=0}^{m-1}
		\int_{t_i}^{t_{i+1}}
		\left\|
		\Bigl(
		E(t_m-s)
		-
		\mathcal{K}_{m-i}
		\Bigr)
		\bigl(G(u(s))-G(u(t_m))\bigr)
		\right\|_{\mathcal{L}_2^0}^{2}
		\,\mathrm{d}s
		\right)^{1/2}
		\right\|_{L^p(\Omega;\mathbb{R})}
		\nonumber\\
		&+
		C_p\left\|
		\left(
		\sum_{i=0}^{m-1}
		\int_{t_i}^{t_{i+1}}
		\left\|
		\Bigl(
		E(t_m-s)
		-
		\mathcal{K}_{m-i}
		\Bigr)G(u(t_m))
		\right\|_{\mathcal{L}_2^0}^{2}
		\,\mathrm{d}s
		\right)^{1/2}
		\right\|_{L^p(\Omega;\mathbb{R})}
		\nonumber\\
		& =:
		C_p(\epsilon_{21}+\epsilon_{22}+\epsilon_{23}).
		\label{eq:4.9}
	\end{flalign}
	We next estimate the three terms $\epsilon_{21}$, $\epsilon_{22}$, and $\epsilon_{23}$ separately.
	
	Under the stability condition, the uniform propagator
	bounds in \cite{ABV2023}, together with
	$\|P_h\|_{\mathcal{L}(H)}\leq1$, yield
	$\left\|
	\mathcal{A}_k(-\tau A_h)^{m-i-1}
	\mathcal{B}_k(-\tau A_h)P_h
	\right\|_{\mathcal{L}(H)}
	\leq C,~i=0,\ldots,m-1, $
	where $C$ is independent of $h$, $\tau$, $k$, and $m$.
	Combining this bound with Assumption~\ref{ass:noise},
	Minkowski's inequality, and Theorem~\ref{th1},
	we obtain
	\begin{flalign}
		\epsilon_{21}
		&\leq
		C\left\|
		\left(
		\sum_{i=0}^{m-1}
		\int_{t_i}^{t_{i+1}}
		\|u(s)-u_{i,h}\|^2\,\mathrm{d}s
		\right)^{1/2}
		\right\|_{L^p(\Omega;\mathbb{R})}
		\nonumber\\
		&\leq
		C\left(
		\sum_{i=0}^{m-1}
		\int_{t_i}^{t_{i+1}}
		\|u(s)-u_{i,h}\|_{L^p(\Omega;H)}^2
		\,\mathrm{d}s
		\right)^{1/2}
		\nonumber\\
		&\leq
		C\left(
		\tau\sum_{i=0}^{m-1}
		\|u(t_i)-u_{i,h}\|_{L^p(\Omega;H)}^2
		\right)^{1/2}
		+
		C\left(
		\sum_{i=0}^{m-1}
		\int_{t_i}^{t_{i+1}}
		\|u(s)-u(t_i)\|_{L^p(\Omega;H)}^2
		\,\mathrm{d}s
		\right)^{1/2}
		\nonumber\\
		&\leq
		C\left(
		\tau\sum_{i=0}^{m-1}
		\|u(t_i)-u_{i,h}\|_{L^p(\Omega;H)}^2
		\right)^{1/2}
		+
		C\left(
		\sum_{i=0}^{m-1}
		\int_{t_i}^{t_{i+1}}
		(s-t_i)\,\mathrm{d}s
		\right)^{1/2}
		\nonumber\\
		&\leq
		C\tau^{1/2}
		+
		C\left(
		\tau\sum_{i=0}^{m-1}
		\|u(t_i)-u_{i,h}\|_{L^p(\Omega;H)}^2
		\right)^{1/2}.
		\label{eq:4.10}
	\end{flalign}
	
	For \(\epsilon_{22}\), interpolating Lemma~\ref{lemma3.6zx} between
	\(\theta=0\) and \(\theta=1\) at \(\theta=1-r\), together with
	Assumption~\ref{ass:noise}, Minkowski's inequality, and
	Theorem~\ref{th1}, we obtain
	\begin{flalign}
		\epsilon_{22}
		&\leq
		C\left(h^{1+r}+\tau^{(1+r)/2}\right)
		\left(
		\sum_{i=0}^{m-1}
		\int_{t_i}^{t_{i+1}}
		(t_m-s)^{-1-r}
		\|G(u(s))-G(u(t_m))\|_
		{L^p(\Omega;\mathcal{L}_2^0)}^2
		\,\mathrm{d}s
		\right)^{1/2}
		\nonumber\\
		&\leq
		C\left(h^{1+r}+\tau^{(1+r)/2}\right)
		\left(
		\int_0^{t_m}
		(t_m-s)^{-1-r}
		\|u(s)-u(t_m)\|_{L^p(\Omega;H)}^2
		\,\mathrm{d}s
		\right)^{1/2}
		\nonumber\\
		&\leq
		C\left(h^{1+r}+\tau^{(1+r)/2}\right)
		\left(
		\int_0^{t_m}(t_m-s)^{-r}\,\mathrm{d}s
		\right)^{1/2}
		\nonumber\\
		&=
		C\left(h^{1+r}+\tau^{(1+r)/2}\right)
		\left(\frac{t_m^{1-r}}{1-r}\right)^{1/2}
		\nonumber\\
		&\leq
		C\left(h^{1+r}+\tau^{(1+r)/2}\right),
		\label{eq:4.11}
	\end{flalign}
	where the second inequality uses the weaker bound
	\((t_m-s)^{-1}\le(t_m-s)^{-1-r}\), valid since \(t_m-s\le1\).
	
	For $\epsilon_{23}$, let $(\psi_\ell)_{\ell\geq1}$
	be an orthonormal basis of $H_0$.
	Applying \eqref{eq:3.4.33} with $\beta=1+r$,
	together with Assumption~\ref{ass:noise}
	and Theorem~\ref{th1}, we obtain
	\begin{flalign}
		\epsilon_{23}
		&=
		\left\|
		\left(
		\sum_{\ell=1}^{\infty}
		\sum_{j=1}^{m}
		\int_{t_{j-1}}^{t_j}
		\left\|
		\Bigl(
		E(s)
		-
		\mathcal{K}_{j}
		\Bigr)G(u(t_m))\psi_\ell
		\right\|^2
		\,\mathrm{d}s
		\right)^{1/2}
		\right\|_{L^p(\Omega;\mathbb{R})}
		\nonumber\\
		&\leq
		C\left(h^{1+r}+\tau^{(1+r)/2}\right)
		\left\|
		\left(
		\sum_{\ell=1}^{\infty}
		\|G(u(t_m))\psi_\ell\|_r^2
		\right)^{1/2}
		\right\|_{L^p(\Omega;\mathbb{R})}
		\nonumber\\
		&=
		C\left(h^{1+r}+\tau^{(1+r)/2}\right)
		\|A^{r/2}G(u(t_m))\|_
		{L^p(\Omega;\mathcal{L}_2^0)}
		\nonumber\\
		&\leq
		C\left(h^{1+r}+\tau^{(1+r)/2}\right)
		\left(
		1+\|u(t_m)\|_{L^p(\Omega;\dot{H}^{r})}
		\right)
		\nonumber\\
		&\leq
		C\left(h^{1+r}+\tau^{(1+r)/2}\right).
		\label{eq:4.12}
	\end{flalign}
	
	Set $\varphi_m=\|u(t_m)-u_{h,m}\|_{L^p(\Omega;H)}^2$
	for $m=0,1,\ldots,M$.
	Combining \eqref{eq:4.3}--\eqref{eq:4.12}, we obtain
	\begin{flalign}
		\varphi_m
		\leq
		C\left(
		h^{1+r}+\tau^{1/2}+\tau^{(1+r)/2}
		\right)^2+
		C\left(
		\tau\sum_{i=0}^{m-1}
		t_{m-i}^{-(1-r)/2}\varphi_i^{1/2}
		\right)^2
		+
		C\tau\sum_{i=0}^{m-1}\varphi_i.
		\label{eq:4.13}
	\end{flalign}
	By H\"older's inequality, we have
	\begin{flalign*}
		\left(
		\tau\sum_{i=0}^{m-1}
		t_{m-i}^{-(1-r)/2}\varphi_i^{1/2}
		\right)^2
		&\leq
		\left(
		\tau\sum_{i=0}^{m-1}
		t_{m-i}^{-(1-r)/2}
		\right)
		\left(
		\tau\sum_{i=0}^{m-1}
		t_{m-i}^{-(1-r)/2}\varphi_i
		\right)
		\\
		&\leq
		\frac{2}{1+r}T^{(1+r)/2}
		\tau\sum_{i=0}^{m-1}
		t_{m-i}^{-(1-r)/2}\varphi_i,
	\end{flalign*}
	where we used
	$
	\tau\sum_{j=1}^{m}t_j^{-(1-r)/2}
	\leq
	\int_0^{t_m}s^{-(1-r)/2}\,\mathrm{d}s
	\leq
	\frac{2}{1+r}T^{(1+r)/2}.
	$
	Moreover, since $t_{m-i}\leq T$,
	$$
	\tau\sum_{i=0}^{m-1}\varphi_i
	\leq
	T^{(1-r)/2}
	\tau\sum_{i=0}^{m-1}
	t_{m-i}^{-(1-r)/2}\varphi_i.
	$$
	
	Since $u_{0,h}=P_hu_0$, the standard finite element
	projection estimate gives
	$$
	\varphi_0
	=
	\|(I-P_h)u_0\|_{L^p(\Omega;H)}^2
	\leq
	Ch^{2+2r}
	\|u_0\|_{L^p(\Omega;\dot{H}^{1+r})}^2
	\leq Ch^{2+2r}.
	$$
	Consequently, using $t_m\geq\tau$ and $\tau\leq T$,
	$$
	\tau t_m^{-(1-r)/2}\varphi_0
	\leq
	\tau^{(1+r)/2}\varphi_0
	\leq Ch^{2+2r}.
	$$
	Combining these bounds with \eqref{eq:4.13}
	and using $\tau^{1+r}\leq T^r\tau$, we obtain
	\begin{equation}\label{eq:4.14}
		\varphi_m
		\leq
		C\left(h^{2+2r}+\tau\right)
		+
		C\tau\sum_{i=1}^{m-1}
		t_{m-i}^{-(1-r)/2}\varphi_i,
		\qquad m=1,2,\ldots,M.
	\end{equation}
	Applying Lemma~\ref{lem:4.1} with $\beta=(1+r)/2$
	yields
	$$
	\varphi_m\leq C\left(h^{2+2r}+\tau\right),
	\qquad m=1,2,\ldots,M.
	$$
	Taking square roots yields the desired estimate
	for $m=1,2,\ldots,M$.
	The case $m=0$ follows directly from
	$u_{0,h}=P_hu_0$ and the initial projection estimate.
	Hence,
	$$
	\|u(t_m)-u_{h,m}\|_{L^p(\Omega;H)}
	\leq C\left(h^{1+r}+\tau^{1/2}\right),
	\qquad m=0,1,\ldots,M.
	$$
	The constant $C$ is independent of $h$, $\tau$, $m$,
	and the number of internal stages $k$.
	This completes the proof.
\end{proof}

\section{Improved strong convergence rates for additive noise} \label{sec5}

For spatially regular additive noise and \(F=0\), \cite{ABV2023}
established a strong convergence rate of order \(1-\epsilon\) in time.
In this section, we develop a complete analysis
for semilinear problems under additional regularity
and differentiability assumptions.
Using the propagator estimates established in
Section~3, we obtain the improved strong error bound
$\mathcal{O}(h^{1+r}+\tau^{(1+r)/2})$
for $r\in(-1,1]$.
In particular, $r=1$ yields second-order spatial
and first-order temporal strong convergence
without an $\varepsilon$-loss.

In this section, we consider additive noise with a
deterministic coefficient $\Phi(t)$ that is independent
of the solution. The equation takes the form \cite{DaPratoZabczyk2014}
\begin{equation}\label{eq:additive-model}
	\begin{cases}
		\mathrm{d}u(t)+Au(t)\,\mathrm{d}t
		=
		F(u(t))\,\mathrm{d}t
		+
		\Phi(t)\,\mathrm{d}W^Q(t),
		& t\in(0,T],\\
		u(0)=u_0.
	\end{cases}
\end{equation}
The time-independent case is included by taking
$\Phi(t)\equiv G_0$.
The assumptions on $F$, $\Phi$, and $u_0$
used in this section are stated below.

The assumptions on $A$, the stochastic basis, and
the spatial discretization remain unchanged.
For the additive-noise analysis, we replace the
conditions on the drift, the diffusion coefficient,
and the initial value by the following assumptions.
Throughout this section, let $r\in(-1,1]$,
so that the spatial regularity index $1+r$
belongs to $(0,2]$.

\begin{assumption}
	\label{ass:additive-F}
	The nonlinear operator $F:H\to H$, there exist constants $L>0$ and $\eta\in[1,2)$
	such that
	\begin{align}
		\|F(u)\|
		&\leq L(1+\|u\|),
		&& u\in H,
		\label{eq:additive-F-growth}
		\\
		\|F'(u)v\|
		&\leq L\|v\|,
		&& u,v\in H,
		\label{eq:additive-F-first}
		\\
		\|F''(u)(v_1,v_2)\|_{-\eta}
		&\leq L\|v_1\|\|v_2\|,
		&& u,v_1,v_2\in H.
		\label{eq:additive-F-second}
	\end{align}
\end{assumption}

In contrast to Assumption~\ref{non}, the drift takes values in $H$ and satisfies additional
differentiability conditions. These additional assumptions are used to obtain
improved temporal convergence rates for additive noise. For further details and examples of nonlinearities satisfying these conditions, we refer to  \cite{Wang2017}.

\begin{assumption}
	\label{ass:additive-noise}
	The noise coefficient is independent of the
	solution and is given by a deterministic mapping
	$ 	\Phi:[0,T]\to	\mathcal{L}_2(H_0,\dot{H}^{r}).$
	There exist constants $L_\Phi>0$ and
	$\delta\in[(1+r)/2,1]$ such that
	\begin{align}
		\|A^{r/2}\Phi(t)\|_{\mathcal{L}_2^0}
		&\leq L_\Phi,
		&& t\in[0,T],
		\label{eq:additive-noise-space}
		\\
		\|A^{r/2}
		(\Phi(t)-\Phi(s))\|_{\mathcal{L}_2^0}
		&\leq L_\Phi|t-s|^\delta,
		&& s,t\in[0,T].
		\label{eq:additive-noise-time}
	\end{align}
\end{assumption}

The time-independent case is included by taking
$\Phi(t)\equiv\Phi_0$, in which case
\eqref{eq:additive-noise-time} holds automatically.
For $r<0$, the coefficient $\Phi(t)$ need not
itself belong to $\mathcal{L}_2^0$;
the spatial regularity condition is imposed
on the weighted operator $A^{r/2}\Phi(t)$.

\begin{assumption}
	\label{ass:additive-init}
	The initial value
	$u_0:\Omega\to\dot{H}^{1+r}$
	is $\mathcal{F}_0/\mathcal{B}(\dot{H}^{1+r})$-
	measurable and satisfies
	$ u_0\in L^p(\Omega;\dot{H}^{1+r})$	for $p\in[2,\infty).$
\end{assumption}

We now state the improved strong convergence result
for additive noise.
\begin{theorem}[Improved strong convergence rates]
	\label{th:additive}
	Let $r\in(-1,1]$ and $p\in[2,\infty)$.
	Suppose that Assumptions~\ref{ass:additive-F},
	\ref{ass:additive-noise}, \ref{ass:additive-init},
	and~\ref{ass:AB} hold.
	Let $k\geq1$ satisfy
	$\tau\lambda_{h,N_h}\leq L_k$.
	Let $u(t)$ be the mild solution of
	\eqref{eq:additive-model}, and let $u_{h,m}$
	be the numerical approximation defined by
	\eqref{eq:3.2.1}, with $G(u_{h,m})$
	replaced by $\Phi(t_m)$.
	Then
	\begin{equation}\label{eq:additive-error}
		\|u(t_m)-u_{h,m}\|_{L^p(\Omega;H)}
		\leq C\left(h^{1+r}+\tau^{(1+r)/2}\right),
		\qquad m=0,1,\ldots,M,
	\end{equation}
	for $h,\tau\in(0,1]$. The constant $C$ is independent of $h$, $\tau$,
	$m$, and the number of internal stages 	satisfying the stability condition.
	In particular, $r=1$ yields second-order spatial and first-order temporal strong convergence
	without an $\varepsilon$-loss.
\end{theorem}

\begin{proof}[Proof of Theorem \ref{th:additive}]
	Set $\zeta=\min\{1+r,1\}$ and write
	$$\mathcal{K}_j = \mathcal{A}_k(-\tau A_h)^{j-1}
	\mathcal{B}_k(-\tau A_h)P_h, \qquad j\geq1.$$
	
	Under the assumptions of this section,
	\cite[Theorem~2.5]{Wang2017} and the growth bound
	on $F$ yield the following estimates for every
	$p\in[2,\infty)$
	\begin{align}
		\sup_{t\in[0,T]}
		\left(
		\|u(t)\|_{L^p(\Omega;\dot{H}^{1+r})}
		+
		\|F(u(t))\|_{L^p(\Omega;H)}
		\right)
		&\leq C,
		\notag\\
		\|u(t)-u(s)\|_{L^p(\Omega;H)}
		&\leq C|t-s|^{\zeta/2},
		\qquad s,t\in[0,T].
		\label{eq:5.1}
	\end{align}
	
	Using the mild solution and the recursive representation
	of the numerical solution, with $u_{0,h}=P_hu_0$, we obtain
	\begin{align}
		\|u(t_m)-u_{h,m}\|_{L^p(\Omega;H)}
		&\leq
		\left\|
		\left(E(t_m)
		-\mathcal{A}_k(-\tau A_h)^mP_h\right)u_0
		\right\|_{L^p(\Omega;H)}
		\notag\\
		&\quad+
		\left\|
		\sum_{i=0}^{m-1}
		\int_{t_i}^{t_{i+1}}
		\left(
		E(t_m-s)F(u(s))
		-\mathcal{K}_{m-i}F(u_{i,h})
		\right)\,\mathrm{d}s
		\right\|_{L^p(\Omega;H)}
		\notag\\
		&\quad+
		\left\|
		\sum_{i=0}^{m-1}
		\int_{t_i}^{t_{i+1}}
		\left(
		E(t_m-s)\Phi(s)
		-\mathcal{K}_{m-i}\Phi(t_i)
		\right)\,\mathrm{d}W^Q(s)
		\right\|_{L^p(\Omega;H)}
		\notag\\
		&:=\epsilon_0+\epsilon_1+\epsilon_2.
		\label{eq:5.2}
	\end{align}
	By \eqref{eq:3.4.4} with $\beta=(1+r)/2$, we have
	\begin{align}
		\epsilon_0
		&\leq
		C\left(h^{1+r}+\tau^{(1+r)/2}\right)
		\|u_0\|_{L^p(\Omega;\dot{H}^{1+r})}
		\notag\\
		&\leq C\left(h^{1+r}+\tau^{(1+r)/2}\right).
		\label{eq:5.3}
	\end{align}
	
	We next decompose the drift error as
	\begin{align}
		\epsilon_1
		&\leq
		\left\|
		\sum_{i=0}^{m-1}
		\int_{t_i}^{t_{i+1}}
		E(t_m-s)\bigl(F(u(s))-F(u(t_i))\bigr)
		\,\mathrm{d}s
		\right\|_{L^p(\Omega;H)}
		\notag\\
		&\quad+
		\left\|
		\sum_{i=0}^{m-1}
		\int_{t_i}^{t_{i+1}}
		\bigl(E(t_m-s)-\mathcal{K}_{m-i}\bigr)
		F(u(t_i))\,\mathrm{d}s
		\right\|_{L^p(\Omega;H)}
		\notag\\
		&\quad+
		\left\|
		\sum_{i=0}^{m-1}
		\int_{t_i}^{t_{i+1}}
		\mathcal{K}_{m-i}
		\bigl(F(u(t_i))-F(u_{i,h})\bigr)
		\,\mathrm{d}s
		\right\|_{L^p(\Omega;H)}
		\notag\\
		&:=\epsilon_{11}+\epsilon_{12}+\epsilon_{13}.
		\label{eq:5.4}
	\end{align}
	
	The uniform boundedness of $\mathcal{K}_j$ and
	Assumption~\ref{ass:additive-F} imply
	\begin{align}
		\epsilon_{13}
		&\leq
		C\sum_{i=0}^{m-1}
		\int_{t_i}^{t_{i+1}}
		\|F(u(t_i))-F(u_{i,h})\|_{L^p(\Omega;H)}
		\,\mathrm{d}s
		\notag\\
		&\leq
		C\tau\sum_{i=0}^{m-1}
		\|u(t_i)-u_{i,h}\|_{L^p(\Omega;H)}.
		\label{eq:5.5}
	\end{align}	
	
	For $\epsilon_{12}$, substituting
	$\sigma=t_m-s$ and $j=m-i$, and using
	Lemma~\ref{lemma3.6zx} with $\theta=0$,
	\eqref{eq:3.4.32} with $\widetilde{\rho}=0$,
	and \eqref{eq:5.1}, we obtain
	\begin{align}
		\epsilon_{12}
		&=
		\left\|
		\sum_{j=1}^{m}
		\int_{t_{j-1}}^{t_j}
		\bigl(E(\sigma)-\mathcal{K}_j\bigr)
		F(u(t_{m-j}))\,\mathrm{d}\sigma
		\right\|_{L^p(\Omega;H)}
		\notag\\
		&\leq
		\left\|
		\sum_{j=2}^{m}
		\int_{t_{j-1}}^{t_j}
		\bigl(E(\sigma)-\mathcal{K}_j\bigr)
		\bigl(F(u(t_{m-j}))-F(u(t_{m-1}))\bigr)
		\,\mathrm{d}\sigma
		\right\|_{L^p(\Omega;H)}
		\notag\\
		&\quad+
		\left\|
		\left(
		\sum_{j=1}^{m}
		\int_{t_{j-1}}^{t_j}
		\bigl(E(\sigma)-\mathcal{K}_j\bigr)
		\,\mathrm{d}\sigma
		\right)F(u(t_{m-1}))
		\right\|_{L^p(\Omega;H)}
		\notag\\
		&\leq
		C(h^2+\tau)
		\sum_{j=2}^{m}
		\int_{t_{j-1}}^{t_j}
		\sigma^{-1}t_{j-1}^{\zeta/2}
		\,\mathrm{d}\sigma
		+C(h^2+\tau)
		\notag\\
		&\leq
		C(h^2+\tau)
		\left(
		\int_0^{t_m}\sigma^{-1+\zeta/2}
		\,\mathrm{d}\sigma+1
		\right)
		\leq C(h^2+\tau).
		\label{eq:5.6}
	\end{align}
	
	For $r>0$, a direct use of the Lipschitz continuity
	of $F$ and the temporal regularity of $u$ yields
	only an $\mathcal{O}(\tau^{1/2})$ bound,
	which is insufficient for the desired convergence rate.
	We therefore follow established arguments based on
	Taylor expansions , as used in \cite{Wang2017}.
	To estimate $\epsilon_{11}$, we use Taylor's formula
	and the mild representation on $[t_i,s]$ to write
	\begin{align}
		F(u(s))-F(u(t_i))
		&=
		F'(u(t_i))(E(s-t_i)-I)u(t_i)
		\notag\\
		&\quad+
		F'(u(t_i))
		\int_{t_i}^{s}E(s-\sigma)F(u(\sigma))
		\,\mathrm{d}\sigma
		\notag\\
		&\quad+
		F'(u(t_i))
		\int_{t_i}^{s}E(s-\sigma)\Phi(\sigma)
		\,\mathrm{d}W^Q(\sigma)
		+R_i(s),
		\label{eq:5.7}
	\end{align}
	where Assumption~\ref{ass:additive-F} and
	\eqref{eq:5.1} yield
	\begin{align}
		\|R_i(s)\|_{L^p(\Omega;\dot{H}^{-\eta})}
		&\leq
		C\|u(s)-u(t_i)\|_{L^{2p}(\Omega;H)}^2
		\notag\\
		&\leq C(s-t_i)^\zeta.
		\label{eq:5.8}
	\end{align}
	Set
	$$
	Z_i=
	\int_{t_i}^{t_{i+1}}
	E(t_m-s)F'(u(t_i))
	\left(
	\int_{t_i}^{s}
	E(s-\sigma)\Phi(\sigma)
	\,\mathrm{d}W^Q(\sigma)
	\right)\,\mathrm{d}s.
	$$
	Using semigroup smoothing and \eqref{eq:5.1},
	we obtain from \eqref{eq:5.7}--\eqref{eq:5.8}
	\begin{align}
		\epsilon_{11}
		&\leq
		C\sum_{i=0}^{m-1}
		\int_{t_i}^{t_{i+1}}
		\left(
		(s-t_i)^{(1+r)/2}
		+(s-t_i)
		+(t_m-s)^{-\eta/2}(s-t_i)^\zeta
		\right)\,\mathrm{d}s
		\notag\\
		&\quad+
		\left\|\sum_{i=0}^{m-1}Z_i\right\|_{L^p(\Omega;H)}
		\notag\\
		&\leq
		C\left(\tau^{(1+r)/2}+\tau+\tau^\zeta\right)
		+
		\left\|\sum_{i=0}^{m-1}Z_i\right\|_{L^p(\Omega;H)}.
		\label{eq:5.9}
	\end{align}
	Here $\eta<2$ ensures that
	$(t_m-s)^{-\eta/2}$ is integrable.
	Since $Z_i$ is $\mathcal{F}_{t_{i+1}}$-measurable
	and $\mathbb{E}[Z_i\mid\mathcal{F}_{t_i}]=0$,
	the discrete BDG inequality, followed by
	Lemma~\ref{lem:bdg}, gives
	\begin{align}
		\left\|\sum_{i=0}^{m-1}Z_i\right\|_{L^p(\Omega;H)}
		&\leq
		C_p\left(
		\sum_{i=0}^{m-1}
		\|Z_i\|_{L^p(\Omega;H)}^2
		\right)^{1/2}
		\notag\\
		&\leq
		C\left(
		\tau\sum_{i=0}^{m-1}
		\int_{t_i}^{t_{i+1}}\int_{t_i}^{s}
		\|E(s-\sigma)\Phi(\sigma)\|_{\mathcal{L}_2^0}^2
		\,\mathrm{d}\sigma\,\mathrm{d}s
		\right)^{1/2}
		\notag\\
		&\leq
		C\left(
		\tau\sum_{i=0}^{m-1}
		\int_{t_i}^{t_{i+1}}\int_{t_i}^{s}
		(s-\sigma)^{\nu-1}
		\,\mathrm{d}\sigma\,\mathrm{d}s
		\right)^{1/2}
		\notag\\
		&\leq C\tau^{(1+\zeta)/2}.
		\label{eq:5.10}
	\end{align}
	Since $\tau\in(0,1]$ and $r\in(-1,1]$,
	\eqref{eq:5.9}--\eqref{eq:5.10} imply
	\begin{align}
		\epsilon_{11}
		&\leq
		C\left(
		\tau^{(1+r)/2}+\tau+\tau^\zeta
		+\tau^{(1+\zeta)/2}
		\right)
		\leq C\tau^{(1+r)/2}.
		\label{eq:5.11}
	\end{align}
	Combining \eqref{eq:5.5}, \eqref{eq:5.6},
	and \eqref{eq:5.11}, we conclude that
	\begin{align}
		\epsilon_1
		&\leq
		C\tau^{(1+r)/2}+C(h^2+\tau)+
		C\tau\sum_{i=0}^{m-1}
		\|u(t_i)-u_{i,h}\|_{L^p(\Omega;H)}.
		\label{eq:5.12}
	\end{align}
	
	Finally, we estimate the stochastic error.
	Applying Lemma~\ref{lem:bdg} and making the same
	change of variables as above, we obtain
	\begin{align}
		\epsilon_2
		&\leq
		C_p\left(
		\sum_{j=1}^{m}
		\int_{t_{j-1}}^{t_j}
		\|E(\sigma)\Phi(t_m-\sigma)
		-\mathcal{K}_j\Phi(t_{m-j})\|_{\mathcal{L}_2^0}^2
		\,\mathrm{d}\sigma
		\right)^{1/2}
		\notag\\
		&\leq
		C_p\left(
		\sum_{j=1}^{m}
		\int_{t_{j-1}}^{t_j}
		\left\|
		E(\sigma)
		\bigl(\Phi(t_m-\sigma)-\Phi(t_{m-j})\bigr)
		\right\|_{\mathcal{L}_2^0}^2
		\,\mathrm{d}\sigma
		\right)^{1/2}
		\notag\\
		&\quad+
		C_p\left(
		\sum_{j=1}^{m}
		\int_{t_{j-1}}^{t_j}
		\left\|
		\bigl(E(\sigma)-\mathcal{K}_j\bigr)
		\Phi(t_{m-j})
		\right\|_{\mathcal{L}_2^0}^2
		\,\mathrm{d}\sigma
		\right)^{1/2}
		\notag\\
		&:=\epsilon_{21}+\epsilon_{22}.
		\label{eq:5.13}
	\end{align}
	By semigroup smoothing and
	Assumption~\ref{ass:additive-noise}, we have
	\begin{align}
		\epsilon_{21}
		&\leq
		C\left(
		\sum_{j=1}^{m}
		\int_{t_{j-1}}^{t_j}
		\sigma^{\zeta-1}(t_j-\sigma)^{2\delta}
		\,\mathrm{d}\sigma
		\right)^{1/2}
		\notag\\
		&\leq
		C\tau^\delta
		\left(
		\int_0^{t_m}\sigma^{\zeta-1}
		\,\mathrm{d}\sigma
		\right)^{1/2}
		\leq C\tau^{(1+r)/2}.
		\label{eq:5.14}
	\end{align}
	
	For $\epsilon_{22}$, let $(\psi_\ell)_{\ell\geq1}$
	be an orthonormal basis of $H_0$.
	Using \eqref{eq:3.4.19} with $\mu=1+r$,
	\eqref{eq:3.4.33} with $\beta=1+r$,
	and Assumption~\ref{ass:additive-noise}, we obtain
	\begin{align}
		\epsilon_{22}
		&\leq
		C_p\left(
		\sum_{j=2}^{m}
		\int_{t_{j-1}}^{t_j}
		\left\|
		\bigl(E(\sigma)-\mathcal{K}_j\bigr)
		\bigl(\Phi(t_{m-j})-\Phi(t_{m-1})\bigr)
		\right\|_{\mathcal{L}_2^0}^2
		\,\mathrm{d}\sigma
		\right)^{1/2}
		\notag\\
		&\quad+
		C_p\left(
		\sum_{\ell=1}^{\infty}
		\sum_{j=1}^{m}
		\int_{t_{j-1}}^{t_j}
		\left\|
		\bigl(E(\sigma)-\mathcal{K}_j\bigr)
		\Phi(t_{m-1})\psi_\ell
		\right\|^2
		\,\mathrm{d}\sigma
		\right)^{1/2}
		\notag\\
		&\leq
		C\left(h^{1+r}+\tau^{(1+r)/2}\right)
		\left(
		\sum_{j=2}^{m}
		\int_{t_{j-1}}^{t_j}
		\sigma^{-1}
		\left\|
		A^{r/2}
		\bigl(\Phi(t_{m-j})-\Phi(t_{m-1})\bigr)
		\right\|_{\mathcal{L}_2^0}^2
		\,\mathrm{d}\sigma
		\right)^{1/2}
		\notag\\
		&\quad+
		C\left(h^{1+r}+\tau^{(1+r)/2}\right)
		\|A^{r/2}\Phi(t_{m-1})\|_{\mathcal{L}_2^0}
		\notag\\
		&\leq
		C\left(h^{1+r}+\tau^{(1+r)/2}\right)
		\left(
		\left(
		\int_0^{t_m}\sigma^{-1+2\delta}
		\,\mathrm{d}\sigma
		\right)^{1/2}+1
		\right)
		\notag\\
		&\leq
		C\left(h^{1+r}+\tau^{(1+r)/2}\right).
		\label{eq:5.15}
	\end{align}
	
	Combining \eqref{eq:5.2}, \eqref{eq:5.3},
	\eqref{eq:5.12}, \eqref{eq:5.14},
	and \eqref{eq:5.15}, and using
	$h^2+\tau\leq h^{1+r}+\tau^{(1+r)/2}$, we obtain
	\begin{align}
		\|u(t_m)-u_{h,m}\|_{L^p(\Omega;H)}
		&\leq
		C\left(h^{1+r}+\tau^{(1+r)/2}\right)+
		C\tau\sum_{i=0}^{m-1}
		\|u(t_i)-u_{i,h}\|_{L^p(\Omega;H)}.
		\label{eq:5.16}
	\end{align}
	The contribution with $i=0$ satisfies
	$ C\tau\|u_0-P_hu_0\|_{L^p(\Omega;H)} \leq Ch^{1+r} $ and can be absorbed into the first term
	on the right-hand side of \eqref{eq:5.16}.
	Applying Lemma~\ref{lem:4.1} with $\beta=1$
	then yields 
	\begin{equation}\label{eq:5.17}
		\|u(t_m)-u_{h,m}\|_{L^p(\Omega;H)}
		\leq C\left(h^{1+r}+\tau^{(1+r)/2}\right) 		\qquad m=0,1,\ldots,M.
	\end{equation}
	The constant $C$ is independent of $h$, $\tau$, $m$,
	and the admissible number of internal stages $k$.
	This completes the proof.
\end{proof}

		\section{Numerical examples}
		\label{sec6}
		In this section, we consider one-dimensional  and two-dimensional examples
		to verify the strong convergence rates established above and
		to compare the computational efficiency of the explicit
		stabilized methods with that of the linear implicit Euler method. 
		\subsection{Test problems and numerical setup}
		We consider the semilinear stochastic heat equation
		with multiplicative noise on $D=(0,1)^d$, $d=1,2$
		\begin{equation}\label{eq:num-multiplicative}
			\begin{cases}
				\mathrm{d}u(t)
				=
				\bigl[a\Delta u(t)+\sin(u(t))\bigr]
				\,\mathrm{d}t
				+u(t)\,\mathrm{d}W^Q(t),
				& t\in(0,T],\\
				u(t,\boldsymbol{x})=0,
				& \boldsymbol{x}\in\partial D,\\
				u(0,\boldsymbol{x})
				=
				\displaystyle\prod_{\ell=1}^d
				\sin(\pi x_\ell),
				& \boldsymbol{x}\in D.
			\end{cases}
		\end{equation}
		We take $a=0.1$,  and $T=1$. Let
		$$
		e_{\boldsymbol{n}}(\boldsymbol{x})
		=
		2^{d/2}\prod_{\ell=1}^d
		\sin(n_\ell\pi x_\ell),
		\qquad
		\mu_{\boldsymbol{n}}
		=
		\sum_{\ell=1}^d n_\ell^2,
		\qquad \boldsymbol{n}\in\mathbb{N}^d.
		$$
		The covariance operator is diagonal in this basis,
		with
		$$
		Qe_{\boldsymbol{n}}
		=
		q_{\boldsymbol{n}}e_{\boldsymbol{n}},
		\qquad
		q_{\boldsymbol{n}}
		=
		\frac{\mu_{\boldsymbol{n}}^{-(d+1)/2}}
		{[\log(2+\mu_{\boldsymbol{n}})]^2}.
		$$
		The driving Wiener process is represented by
		$$
		W^Q(t)
		=
		\sum_{\boldsymbol{n}\in\mathbb{N}^d}
		\sqrt{q_{\boldsymbol{n}}}\,
		\beta_{\boldsymbol{n}}(t)e_{\boldsymbol{n}},
		$$
		where the $\beta_{\boldsymbol{n}}$ are independent
		standard real-valued Brownian motions. It is straightforward to verify that the drift
		$F(v)=\sin(v)$, the diffusion coefficient $G(v)\xi=v\xi$, and the initial value satisfy
		Assumptions~\ref{non}, \ref{ass:noise}, and~\ref{ass:init} with $r=1/2$.
		
		We discretize the Laplacian using standard central finite
		differences\footnote{The standard finite difference scheme
			\eqref{eq:6.2}--\eqref{eq:6.3} on the uniform mesh is the simplest
			\(P_1\) finite element method, with lumped mass matrix \(M_h=h^dI\);
			its operator \(A_h=M_h^{-1}K_h\) is exactly the \((2d+1)\)-point
			stencil, so the analysis of Sections~\ref{sec2}--\ref{sec3} applies.}
		on uniform grids with mesh width \(h=1/N\). In one dimension, writing
		$u_i(t)\approx u(t,x_i)$ with $x_i=ih$, we obtain
		\begin{equation}\label{eq:6.2}
			\begin{aligned}
				\mathrm{d}u_i(t)
				=
				\left[
				a\frac{u_{i+1}(t)-2u_i(t)+u_{i-1}(t)}{h^2}
				+\sin(u_i(t))
				\right]\mathrm{d}t 
				+
				u_i(t)\,\mathrm{d}W_J^Q(t,x_i),
			\end{aligned}
		\end{equation}
		where $i=1,\ldots,N-1$ and $W_J^Q$ is the
		truncated Wiener process defined above,
		with independent standard Wiener processes
		$\beta_n$, $n=1,\ldots,J$.
		The boundary and initial conditions are
		$u_0(t)=u_N(t)=0$ and
		$u_i(0)=\sin(\pi x_i)$.
		
		In two dimensions, writing
		$u_{i,j}(t)\approx u(t,x_i,y_j)$ with
		$x_i=ih$ and $y_j=jh$, we obtain
		\begin{equation}\label{eq:6.3}
			\begin{aligned}
				\mathrm{d}u_{i,j}(t)
				&=
				\Biggl[
				\frac{a}{h^2}
				\Bigl(
				u_{i+1,j}(t)+u_{i-1,j}(t)
				+u_{i,j+1}(t)+u_{i,j-1}(t)
				-4u_{i,j}(t)
				\Bigr)
				\\
				&\qquad+
				\sin(u_{i,j}(t))
				\Biggr]\mathrm{d}t
				+
				u_{i,j}(t)\,
				\mathrm{d}W_J^Q(t,x_i,y_j),
			\end{aligned}
		\end{equation}
		where $i,j=1,\ldots,N-1$, the boundary and initial conditions are
		$u_0(t)=u_N(t)=0$ and $u_{i,j}(0)=\sin(\pi x_i)\sin(\pi y_j)$.
		For time integration, we use the SK-ROCK method \eqref{eq:3.2.9}
		and its variants  \eqref{eq:3.2.12} described in Section~\ref{sec3}.
		The linearly implicit Euler method is also
		included for comparison.
		For SK-ROCK and the variants considered here,
		we take the damping parameter $\eta=0.05$, an initial estimate of the required number
		of internal stages is \cite{AbdulleAlmuslimaniVilmart2018}
		\begin{equation}\label{eq:6.4}
			k=
			\left\lceil
			\sqrt{
				\frac{\tau\rho_h+1.5}
				{2-\frac{4}{3}\eta}
			}
			\right\rceil,
			\qquad
			\rho_h=4ad h^{-2},
		\end{equation}
		where $\rho_h$ is an upper bound for the
		spectral radius of the discrete diffusion
		operator.

		\subsection{Strong convergence tests and computational efficiency}
\begin{figure}[htbp]
			\centering
			\begin{minipage}[b]{0.49\textwidth}
				\includegraphics[width=\linewidth]{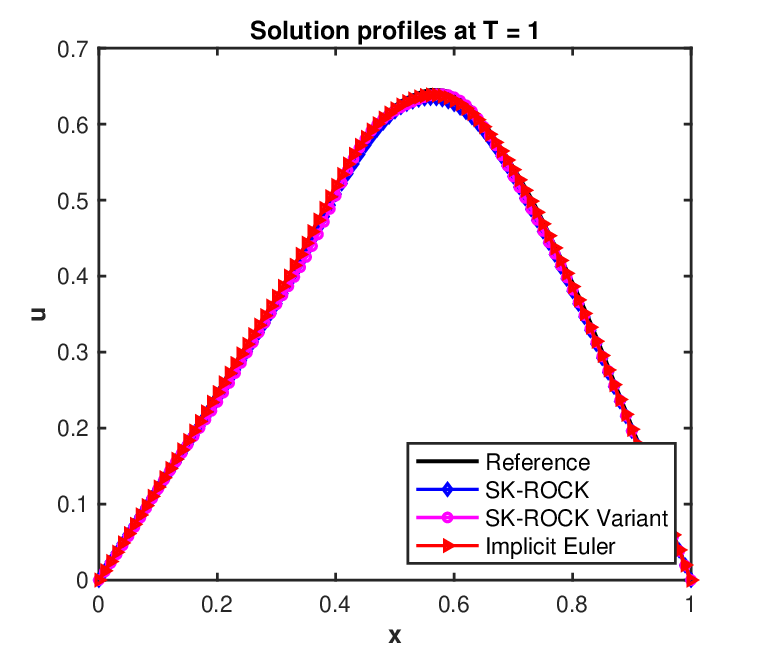} 		
			\end{minipage}
			\hfill 
			\begin{minipage}[b]{0.49\textwidth}
				\includegraphics[width=\linewidth]{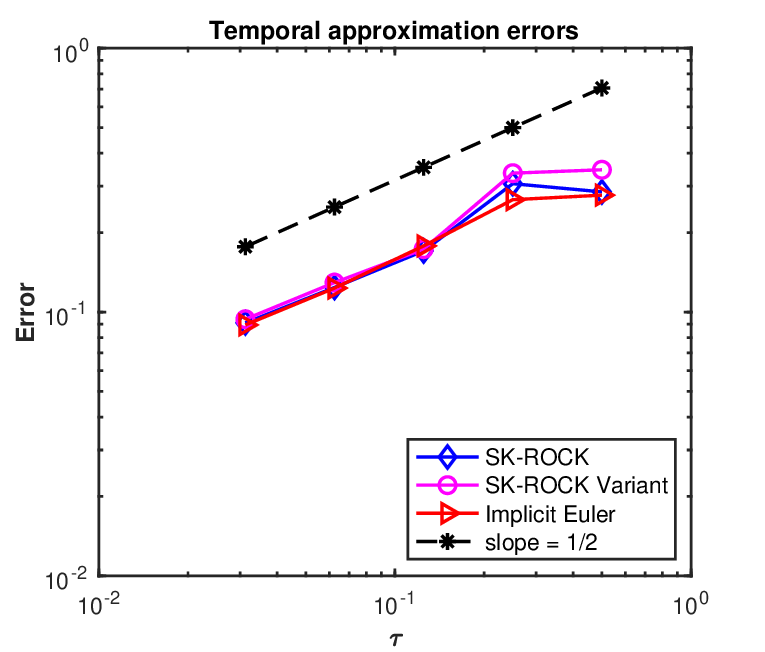} 		 
			\end{minipage}
			\caption{Solution profiles and temporal strong convergence
				for the one-dimensional stochastic heat equation with multiplicative
				noise: spatial profiles at $T=1$ (left) and strong errors versus
				time step size $\tau$ (right).}
			\label{fig:ex1}
		\end{figure}

		We now examine the strong convergence and computational efficiency
		of the three methods.  Fig.~\ref{fig:ex1} presents the numerical results
		for the one-dimensional problem with multiplicative noise.
		We fix the spatial mesh width at $h=1/100$ and use
		the time steps $\tau=2^{-j}$, $j=1,\ldots,5$,
		for the temporal convergence tests.
		The reference solution is computed on the same spatial
		grid using the linearly implicit Euler method with
		$\tau_{\mathrm{ref}}=2^{-15}$.
		All approximations are driven by the same Wiener paths,
		with coarse-grid increments obtained by summing
		the corresponding fine-grid increments.
		The left panel compares the solution profiles at $T=1$
		for a single noise realization and $\tau=2^{-5}$.
		All three numerical solutions closely follow
		the reference profile.
		The right panel displays the root-mean-square errors
		in the discrete $L^2$-norm, estimated from $100$
		independent samples.
		For the smaller time steps, the error curves are
		approximately parallel to the reference line of
		slope $1/2$, consistent with the theoretical
		temporal strong convergence rate. Since the efficiency gains of explicit stabilized methods
		can be more pronounced in higher spatial dimensions,
		we compare the computational costs of the three methods
		in the two-dimensional example below.
		
\begin{figure}[htbp]
			\centering
			\begin{minipage}[b]{0.49\textwidth}
				\includegraphics[width=\linewidth]{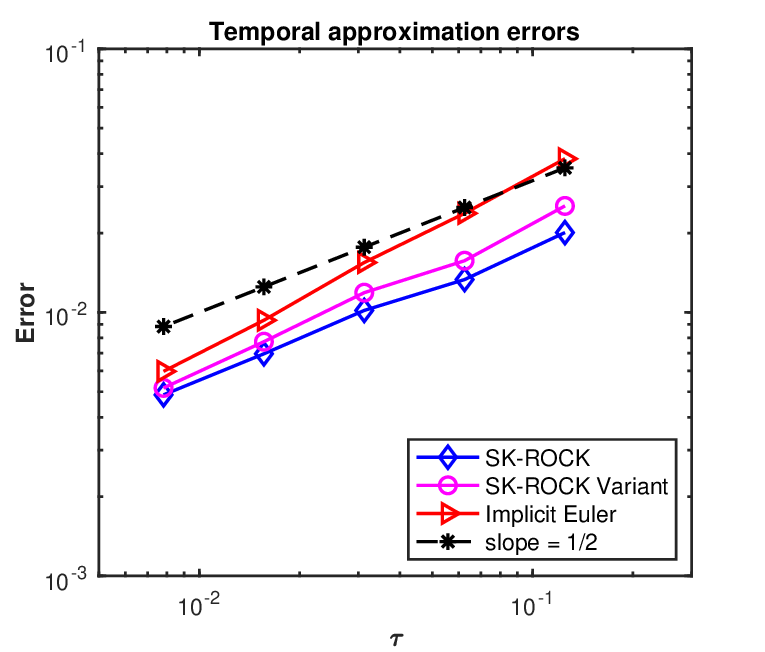} 		
			\end{minipage}
			\hfill 
			\begin{minipage}[b]{0.49\textwidth}
				\includegraphics[width=\linewidth]{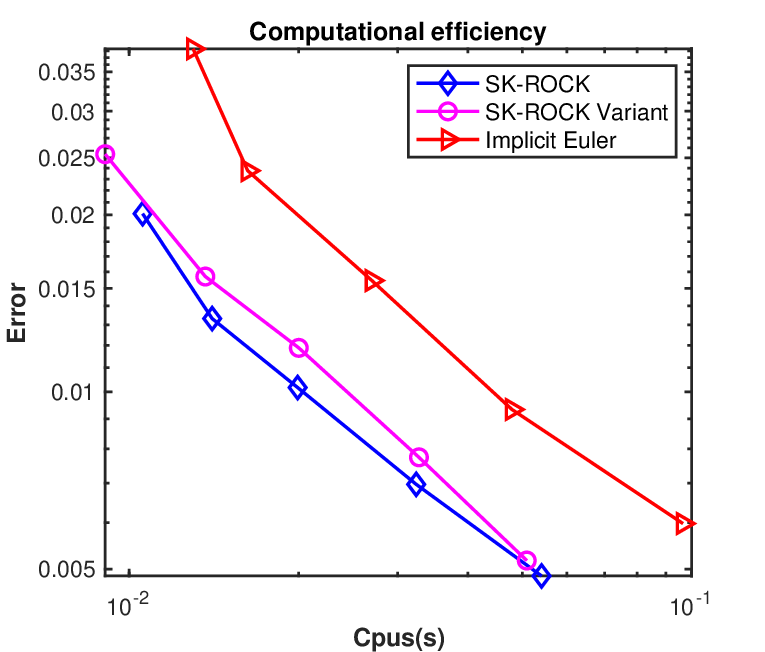} 		 
			\end{minipage}
			\caption{Temporal strong convergence and computational efficiency
				for the two-dimensional stochastic heat equation with multiplicative
				noise: strong errors versus time step size (left) and CPU time (right).}
			\label{fig:ex2}
		\end{figure}

		We next consider the two-dimensional problem with
		multiplicative noise.
		Fig.~\ref{fig:ex2} presents the temporal strong errors
		and computational efficiency of the three methods.
		We fix $h=1/100$ and use the time steps
		$\tau=2^{-j}$, $j=3,\ldots,7$.
		The reference solution is computed on the same spatial
		grid with $\tau_{\mathrm{ref}}=2^{-15}$,
		and all approximations are driven by the same
		coupled Wiener paths.
		The left panel shows the root-mean-square errors
		at $T=1$ in the discrete $L^2$-norm.
		The error curves approximately follow the reference
		slope $1/2$ for the smaller time steps,
		consistent with the theoretical temporal strong
		convergence rate.
		The right panel compares the errors against CPU time.
		Both explicit stabilized methods achieve comparable
		accuracy at a lower computational cost than the
		linearly implicit Euler method, with SK-ROCK generally
		providing the most favorable balance between
		accuracy and computational cost in this experiment. 

		\section{Conclusion}
		
		We have derived estimates for the discrete propagation
		operators that are uniform with respect to the  number of internal stages $k$. Based on these estimates, we established the strong error
		bound $\mathcal{O}(h^{1+r}+\tau^{1/2})$ for fully discrete
		explicit stabilized methods applied to semilinear parabolic
		SPDEs with multiplicative noise.
		Under additional regularity assumptions, we obtain convergence orders without $\epsilon$-loss for the special case of additive noise.
		Future work will extend the analysis to more general nonlinear
		SPDEs, including stochastic Allen--Cahn equations,
		and to equations driven by fractional Brownian motion.
		
		\section*{Acknowledgements}
		The computations in this paper were performed in part on the High Performance Computing
		Platform of Xiangtan University.
		
\section*{Funding}
This work was supported by the National Natural Science Foundation of China
(No.~12471391, to A.~Xiao), the Natural Science Foundation of Hunan Province
(No.~2023JJ40615, to X.~Tang), and the Excellent Youth Foundation of the
Education Department of Hunan Province (No.~25B0151, to X.~Tang).
		
		\bibliographystyle{abbrvnat}  
		\bibliography{refs2}
	\end{document}